\documentclass[12pt]{amsart}
\usepackage[utf8]{inputenc}
\usepackage{amssymb}
\usepackage[active]{srcltx}       
\numberwithin{equation}{section}

\title{List chromatic number of finitary matroids: A generalization of Seymour's result}
\thanks{Supported  by  NKFIH grant K129211. } 
\author[T. Csernák]{Tamás Csernák} 
\address{Eötvös University of Budapest, Hungary}
\email{tamas@csernak.com}

\usepackage{amsthm}
\usepackage{amsfonts}

\usepackage{mathrsfs}
\usepackage{enumerate}
\makeatletter
\@namedef{subjclassname@2020}{%
  \textup{2020} Mathematics Subject Classification}
\makeatother
\newtheorem{theorem}{Theorem}[section]        

\newtheorem{lemma}[theorem]{Lemma}
\newtheorem{observation}[theorem]{Observation}
\newtheorem{corollary}[theorem]{Corollary}
\theoremstyle{definition}
\newtheorem{definition}[theorem]{Definition} 

\newtheorem{problem}[theorem]{Problem}

\newtheorem{example}[theorem]{Example}

\theoremstyle{remark}
\newtheorem{remark}[theorem]{Remark}

\newcommand{\bigslant}[2]{{\raisebox{.2em}{$#1$}\left/\raisebox{-.2em}{$#2$}\right.}}

\newcommand{\Dom}{\operatorname{Dom}}

\begin{document}

\begin{abstract}
 Seymour proved that the 
 chromatic numbers and the list chromatic numbers of any loop-free finite matroids 
 are the same. In this paper we prove the same statement for infinite, 
 loop-free finitary matroids. 
\end{abstract}

\maketitle

\section{Introduction}

Matroids are important objects in finite combinatorics that can represent the basic properties of independence and rank. Theorems related to  matroids have  application in many areas of mathematics, 
such as linear algebra or graph theory. One of the most interesting fact about matroids is Seymour's list coloring theorem \cite{Se98}, which  states that 
 the list chromatic number of a finite matroid is  equal to the chromatic number.  
 
 In this paper, we first generalize this result to a class of infinite matroids called finitary matroid. In section \ref{MAT_sc:def} we describe the most important properties of these matroids. Then, in  section \ref{MAT_sc:fin}, we prove the generalization of Seymour's theorem, when the chromatic number is finite, by using some more or less standard compactness arguments. In section \ref{MAT_sc:inf}, we show the same statement when the chromatic number is an infinite cardinal. This proof uses  heavier set theory and logic, in particular,  elementary submodels. 

 Another result by Bartinski et all. in \cite{Coop} gave a generalization of Seymour's result in the case where multiple matroids were involved.

\begin{theorem}\label{TH_CooperativeFinite}
Theorem 2 in \cite{Coop}
    If $\mathcal{M}_{1},...,\mathcal{M}_{t}$ are loop-free finite matroids on the same base set $S$ and there is some $k$, such that $Chr(\mathcal{M}_{i})\le k$ and a listing $L$ from $S$, with $L(x)\subset \{1,...,t\}, |L(x)|\ge k$ for all $x\in S$. Then there is a $\Phi: S\rightarrow \{1,...,t\}$ function, such that

    1. $\Phi (x)\in L(x)$ for all $x\in S$

    2. $\Phi ^{-1}(i)\in \mathcal{I}(\mathcal{M}_{i})$ for all $1\le i\le t$.
\end{theorem}

In section \ref{MAT:sc_coop} we provide this generalization for finitary matroids. Later, in section \ref{MAT_sc:red} we initiate a strong tool called {\em partition reduction} and finally in section \ref{MAT_sc:dual} we will also show some results of Seymour's theorem for the duals of finitary matroids, that are not nessecerly finitary.

For finitary matroids, we use the basic notions and results of \cite{Ox92} for finite matroids. For more general infinite matroids, we will use the results from \cite{Br13}.

\section{Definition and basic properties of finitary matroids}\label{MAT_sc:def}
 
Finitary matroids can be defined various ways, such as by rank function, independent sets, circuits, bases or closure operator. All of these definitions can be shown to be equvivalent. In this paper, the most convenient way is to define matroids by the rank function. We follow the terminology of \cite{Ox92}.  

\begin{definition}\label{MAT_df:finitary-matroid}
 Let $S$ be any set. A {\em rank function on $S$} is a function $r:[S]^{< \omega }\rightarrow \omega $, with the following properties:
\begin{enumerate}[1.]
\item 
 $r(\emptyset )=0$, 
\item  
 $\forall A,B\in [S]^{<\omega }$ if $A\subseteq B$, then $r(A)\le r(B)$ ({\em monotonity}),
\item
  $\forall A\in [S]^{<\omega }$, $r(A)\le |A|$ ({\em subcardinality}),
\item 
 $\forall A,B\in [S]^{<\omega }$, $r(A)+r(B)\ge r(A\cap B)+r(A\cup B)$ ({\em submodularity}).
\end{enumerate}
 
 A {\em finitary matroid} is a pair $\mathcal{M}=(S,r)$, where $r$ is a rank function on $S$. A subset $X\subseteq S$ is called {\em independent} if $\forall A\in [X]^{<\omega }$, we have $r(A)=|A|$. The set of independent sets in $\mathcal{M}$ is denoted by $\mathcal{I}(\mathcal{M})$.
\end{definition}

By this definition it is clear that the subsets of an independent set are also independent.

If $\mathcal{M}=(S,r)$ is a finitary matroid and $A\in [S]^{<\omega }$, then $\mathcal{M}_{A}=(A,r|_{\mathcal{P}(A)})$ is a finite matroid. Using this observation, there are some claims, that can be proven easily using the fact that these hold for finite matroids. On the other hand, if $r$ is any function on the finite subsets of $S$, such that all $\mathcal{M}_{A}$-s are finite matroids, then $\mathcal{M}$ is a finitary matroid. This way, we can define several matroids. For example let $V$ be any vector space and $S\subseteq V$. The {\em linear matroid} of $S$ is a matroid, where each $A\in [S]^{<\omega }$ is independent if and only if their elements are linearly independent. Then $r(A)$ is the dimension of the subspace, generated by the elements of $A$. The other example is the {\em graphical matroid}, where $V$ is a vertex set and $S\subseteq [V]^{2}$ is an edge set. In this matroid   $A\in [S]^{<\omega }$ is independent, if it does not contain a circuit and $r(A)$ is the number of vertices covered by $A$ minus the number of components of these vertices by the edge set $A$. In fact it can be shown, that graphical matroids are a subclass of linear matroids, although there are some basic terms in matroid theory, which came from graphical matroids.
However, there are many operations of finite matroids, that cannot be used for finitary matroids, such as dualisation. 

Let us remark that there are "natural" infinite matroids defined in a more general way in \cite{Br13} which are not finitary: for example the {\em bond matroid} of an infinite graph cannot be obtained as a finitary matroid, as it may have sets that are not independent, but all of their finite subsets are independent. There are several definitions that are slightly different for finitary matroids and the general case in \cite{Br13}, some of them can be easily seen to be equivalent. However to make it more clear, we will show some properties for finitary matroids using only the definitions given here. Some of these have some form of generalization to general infinite matroids.

\begin{lemma}\label{MAT_lm:lemma1}
 If $\mathcal{M}=(S,r)$ is a finitary matroid and $A\in [S]^{<\omega }$, then $A\in \mathcal{I}(\mathcal{M})$ if and only if $r(A)=|A|$. Moreover, $|B|=r(A)$ for each  maximal independent subset $B\subseteq A$.
\end{lemma}

\begin{proof}
 Use the same result, (Theorem 1.3.2. in \cite{Ox92}) for the finite matroid  $\mathcal{M}_{A}$.
\end{proof}

\begin{definition}
 Let $\mathcal{M}=(S,r)$ be a finitary matroid. A subset $C\subseteq S$ is called a {\em circuit}, if $C\not\in \mathcal{I}(\mathcal{M})$, and $C$ is minimal dependent (i.e. not independent). The set of circuits of $\mathcal{M}$ is denoted by $\mathcal{C}(\mathcal{M})$.
\end{definition}

If $X\subseteq S$ is not independent, then there is some $A\in [X]^{<\omega }$ that is not independent. Then taking the elements one by one, we have that there is a minimal not independent set $C\subseteq A\subseteq X$, and thus $C\in \mathcal{C}(\mathcal{M})$. Hence, all not independent subsets of a finitary matroid contain a circuit and all circuits are finite.

It is also clear that $\emptyset \not\in \mathcal{C}(\mathcal{M})$ and if $C_{1}\neq C_{2}\in \mathcal{C}(\mathcal{M})$, then $C_{1}\not\subseteq C_{2}$. 
The following lemma contains two statements, 
{\em weak and strong circuit elimination axioms}, that hold for infinite finitary matroids as well.



\begin{lemma}\label{MAT_lm:lemma2}

 a)If $\mathcal{M}=(S,r)$ is a finitary matroid, $C_{1}\neq C_{2}\in \mathcal{C}(\mathcal{M})$ and $e\in C_{1}\cap C_{2}$, then there is a $C\in \mathcal{C}(\mathcal{M})$, such that $C\subseteq (C_{1}\cup C_{2})-e$. 
 
 b) If we also have $e_{1}\in C_{1}-C_{2}$, $C$ can be chosen such that $e_{1}\in C$.
\end{lemma}

\begin{proof}
 For a) use Lemma 1.1.3 in \cite{Ox92} for the finite matroid $\mathcal{M}_{C_{1}\cup C_{2}}$. For b) use Proposition 1.4.11 in \cite{Ox92} for the finite matroid $\mathcal{M}_{C_{1}\cup C_{2}}$.
\end{proof}

A circuit consisting of  one element is called a {\em loop}. A finitary matroid is called {\em loop-free} if it does not contain any loop. Equivalently, this means that for all $x\in S$, we have $r(\{ x\})=1$.

\begin{definition}
  Let $\mathcal{M}=(S,r)$ be a finitary matroid. A set $B\in \mathcal{I}(\mathcal{M})$ is a {\em base}, if it is a maximal independent set. The set of bases of $\mathcal{M}$ is denoted by $\mathcal{B}(\mathcal{M})$. 
 
\end{definition}

The existence of bases for finite matroids are clear, as we can add elements to an independent sets one by one, until it gets maximal, however since independence depend just on the finite subsets, by Teichmuller-Tukey lemma, we can show that bases exist in all finitary matroids and all independent sets are contained in a base.
 
\begin{lemma}\label{MAT_lm:lemma3}
 If $\mathcal{M}=(S,r)$ is a finitary matroid, $B\in \mathcal{B}(\mathcal{M})$ and $x\in S-B$, then there is a unique $C\in \mathcal{C}(\mathcal{M})$, such that $C\subseteq B+x$. 
\end{lemma} 

\begin{proof}
 Since $B$ is maximal independent, $B+x$ is not independent, so it must contain a circuit $C$.
 Since $B$ is independent, $x\in C$.
 
 For unicity suppose $C_{1},C_{2}\subseteq B+x$ are circuits. Since no circuit can be the subset of the independent set $B$, we must have $x\in C_{1}\cap C_{2}$. But then by lemma \ref{MAT_lm:lemma2} a), there is a $C\in \mathcal{C}(\mathcal{M})$, such that $C\subseteq (C_{1}\cup C_{2})-x\subseteq B$, that is a contradiction.
\end{proof}

\begin{definition}\label{MAT_df:maincircle}
  Let $\mathcal{M}=(S,r)$ be a finitary matroid, $B\in \mathcal{B}(\mathcal{M})$ and $x\in S-B$. The {\em fundamental circuit} of $x$ on $B$, denoted by $C(B,x)$ is the unique $C\in \mathcal{C}(\mathcal{M})$, such that $C\subseteq B+x$.
\end{definition}

The next notion we need to introduce are {\em closed sets} in finitary matroids. But before, we show some important lemmas. It is clear by submodularity and subcardinality, that if $A\in [S]^{<\omega }$ and $x\in S$, then $r(A+x)$ is either $r(A)$ or $r(A)+1$.

\begin{lemma}\label{MAT_lm:lemma4}
 If $\mathcal{M}=(S,r)$ is a finitary matroid, $A,B\in [S]^{<\omega }$, with $A\subseteq B$, and $x\in S$ is such that $r(A+x)=r(A)$. Then $r(B+x)=r(B)$.
\end{lemma} 

\begin{proof}
 Using the submodularíty for the sets $A+x$ and $B$, we have $r(A+x)+r(B)\ge r(A)+r(B+x)$. By our assumption, then $r(A)+r(B)\ge r(A)+r(B+x)$, so $r(B)\ge r(B+x)$. 
 By the monotonity, we have $r(B)\le r(B+x)$, so $r(B)$ and  $r(B+x)$ are equal.
\end{proof}

\begin{lemma}\label{MAT_lm:lemma5}
  If $\mathcal{M}=(S,r)$ is a finitary matroid, $A\in [S]^{<\omega }$ and $x_{1},...,x_{n}\in S$ such that $r(A+x_{i})=r(A)$ for all $1\le i\le n$, then $r(A\cup \{ x_{1},...,x_{n}\})=r(A)$.
\end{lemma}
 
 \begin{proof}
  Induction on $n$. For $n=1$, it is clear. Suppose this is true for some $n$ and show for $n+1$. Since $A\subseteq A\cup \{ x_{1},...,x_{n}\}$ and $r(A+x_{n+1})=r(A)$, by lemma  \ref{MAT_lm:lemma4}, we have $r(A\cup \{ x_{1},...,x_{n},x_{n+1}\})=r(A\cup \{ x_{1},...,x_{n}\}+x_{n+1})=r(A\cup \{ x_{1},...,x_{n}\})=r(A)$ by the induction hypothesis, so the induction step works.
 \end{proof}
 
For a finite matroid a subset $Z\subseteq S$ is called closed if $r(Z+x)>r(Z)$ for all $x\in S-Z$. For finitary matroids, the rank function if defined just on finite subsets, so we need to refine this definition.

\begin{definition}
  Let $\mathcal{M}=(S,r)$ be a finitary matroid. A subset $Z\subseteq S$ is {\em closed} if for all $Z_{0}\in [Z]^{<\omega }$ and $x\in S-Z$, we have $r(Z_{0}+x)>r(Z_{0})$.
\end{definition}

By lemma  \ref{MAT_lm:lemma4}, this is clearly equivalent to the original definition of closedness for finite matroids. The whole set $S$ is closed for all matroids and $\emptyset $ is closed if and only is $\mathcal{M}$ is loop-free.  

\begin{lemma}\label{MAT_lm:lemma6}
 If $\mathcal{M}=(S,r)$ is a finitary matroid, $I$ is an index set and for all $i\in I$, $Z_{i}\subseteq S$ is a closed subset, then $\bigcap_{i\in I}{Z_{i}}$ is closed.
\end{lemma}
 
 \begin{proof}
 Let $Z_{0}\in [\bigcap_{i\in I}{Z_{i}}]^{<\omega }$ and $x\in S-(\bigcap_{i\in I}{Z_{i}})$. Then there is some $i\in I$, such that $x\not\in Z_{i}$. Since $Z_{i}$ is closed and $Z_{0}\in [Z_{i}]^{<\omega }$, we have $r(Z_{0}+x)>r(Z_{0})$. Since $Z_{0}$ and $x$ was arbitrary, we have that $\bigcap_{i\in I}{Z_{i}}$ is closed.
\end{proof}

\begin{definition}
  Let $\mathcal{M}=(S,r)$ be a finitary matroid and $X\subseteq S$, then the {\em closure of $X$}, denoted by $\sigma (X)$, is defined by 
  $$\sigma (X)=\bigcap \{Z\subseteq S:X\subseteq Z,Z  \text{ is closed}\}.$$
\end{definition}

By lemma \ref{MAT_lm:lemma6}, $\sigma (X)$ is closed and $X\subseteq \sigma (X)$, as it was contained in all members of intersection.

\begin{lemma}\label{MAT_lm:lemma7}
  Let $\mathcal{M}=(S,r)$ be a finitary matroid, then
\begin{enumerate}[a)]
\item 
 For all $X\subseteq S$ if $X\subseteq Z\subseteq S$ and $Z$ is closed, then $\sigma (X)\subseteq Z$,
\item   
 If $X\subseteq Y\subseteq S$, then $\sigma (X)\subseteq \sigma (Y)$,
\item 
 If $X\subseteq S$, the $\sigma (\sigma (X))=\sigma (X)$.
\end{enumerate}  
\end{lemma}

\begin{proof}
  a) By the definition of $\sigma (X)$, $Z$ was an item in the intersection, so $\sigma (X)\subseteq Z$.
  
  b) Since $X\subseteq Y\subseteq \sigma (Y)$ and $\sigma (Y)$ is closed, by a), we have $\sigma (X)\subseteq \sigma (Y)$.
  
  c) On one hand $\sigma (X)\subseteq \sigma (\sigma (X))$ is by definition, on the other hand, since $\sigma (X)\subseteq \sigma (X)$ and $\sigma (X)$ is closed, by a) we must have $\sigma (\sigma (X))\subseteq \sigma (X)$, so they are equal.
\end{proof}

Now we give a characterization of the elements of $\sigma (X)$

\begin{lemma}\label{MAT_lm:lemma8}
  If $\mathcal{M}=(S,r)$ is a finitary matroid, $X\subseteq S$, then 
  $$\sigma (X)=\{ x\in S\ |\ \exists X_{0}\in [X]^{<\omega }\ r(X_{0}+x)=r(X_{0})\}.$$
\end{lemma}

\begin{proof}
  Let $$F=\{ x\in S| \exists X_{0}\in [X]^{<\omega }, r(X_{0}+x)=r(X_{0})\}.$$ 
  First we will show  that $F\subseteq \sigma (X)$. Suppose $x\in F$ and let $X_{0}\in [X]^{<\omega }$, such that $r(X_{0}+x)=r(X_{0})$. Since $X_{0}\subseteq X\subseteq \sigma (X)$, and $\sigma (X)$ is closed, so $x\not\in \sigma (X)$ would mean $r(X_{0}+x)>r(X_{0})$, in contradiction with our assumption. Thus, $x\in \sigma (X)$, so $F\subseteq \sigma (X)$.
  
  Clearly $X\subseteq F$, as for $x\in X$, we can write the singleton $X_{0}=\{x\}$ into the definition of $F$. Now, we need to show that $F$ is closed. Let $Y_{0}\in [F]^{<\omega }$, and list all its elements $Y_{0}=\{ y_{1},...,y_{n}\}$ and let $x\in S-F$. For all $1\le i\le n$, since $y_{i}\in F$, there is an $X_{i}\in [X]^{<\omega }$, such that $r(X_{i}+y_{i})=r(X_{i})$. Since for all $i$, we have $X_{i}\subseteq \bigcup_{j=1}^{n}{X_{j}}$, by lemma  \ref{MAT_lm:lemma4}, we have $r((\bigcup_{j=1}^{n}{X_{j}})+x_{i})=r(\bigcup_{j=1}^{n}{X_{j}})$. Then by lemma  \ref{MAT_lm:lemma5}, we have $r((\bigcup_{j=1}^{n}{X_{j}})\cup Y_{0})=r(\bigcup_{j=1}^{n}{X_{j}})$. Since $x\not\in F$ and $\bigcup_{j=1}^{n}{X_{j}}\in [X]^{<\omega }$, by the definition of $F$, we must have $r((\bigcup_{j=1}^{n}{X_{j}})+x)>r(\bigcup_{j=1}^{n}{X_{j}})$. Then $r((\bigcup_{j=1}^{n}{X_{j}})\cup Y_{0}+x)\ge r((\bigcup_{j=1}^{n}{X_{j}})+x)>r(\bigcup_{j=1}^{n}{X_{j}})=r((\bigcup_{j=1}^{n}{X_{j}})\cup Y_{0})$, so by reverting lemma  \ref{MAT_lm:lemma4}, we get that $r(Y_{0}+x)>r(Y_{0})$. Since $Y_{0}$ and $x$ were arbitrary, $F$ is closed. Then by lemma  \ref{MAT_lm:lemma7}, $\sigma (X)\subseteq F$, and by the first part, ${\sigma}(X)$ and $F$ must be equal.
\end{proof}
  
\begin{lemma}\label{MAT_lm:lemma9}
  If $\mathcal{M}=(S,r)$ is a finitary matroid, $B\subseteq S$, then $B\in \mathcal{B}(\mathcal{M})$ if and only if $B\in \mathcal{I}(\mathcal{M})$ and $\sigma (B)=S$.
\end{lemma}  

\begin{proof}
  If $B\in \mathcal{B}(\mathcal{M})$, then it is clearly independent. Suppose $\sigma (B)\neq S$ and let $x\in S-\sigma (B)$. 
  Then for any $B_{0}\in [B]^{<\omega }$, since $B_{0}\subseteq B\subseteq \sigma (B)$, we have $r(B_{0}+x)>r(B_{0})$, so $r(B_{0}+x)=r(B_{0})+1=|B_{0}|+1=|B_{0}+x|$. 
  Since all finite subsets of $B+x$ are either in this form or subset of $B$, by definition, we have $B+x\in \mathcal{I}(\mathcal{M})$, 
  in contradiction, with the maximality of $B$.
  
  Now suppose $B\in \mathcal{I}(\mathcal{M})$ and $\sigma (B)=S$. We need to show that $B$ is maximal independent. Let $x\in S-B$. Then since $x\in \sigma (B)$, by lemma  \ref{MAT_lm:lemma8}, there is a subset $B_{0}\in [B]^{<\omega }$, such that $r(B_{0}+x)=r(B_{0})=|B_{0}|<|B_{0}+x|$. Since $B_{0}+x\in [B+x]^{<\omega }$, $B+x$ is not independent, so $B$ is maximal.
\end{proof}
  
 In the next step, we define the contraction a set $Z\subseteq S$ in a finitary matroid $\mathcal{M}=(S,r)$. For finite matroids, this is a matroid on the set $S-Z$, defined as for all $A\subseteq S-Z$, $r'(A)=r(A\cup Z)-r(Z)$. However, for finitary matroids, we cannot define the contraction by an infinite set that way, so we need a refined definition.
 
 \begin{definition}
   Let $\mathcal{M}=(S,r)$ be a finitary matroid and $Z\subseteq S$. 
   The {\em contraction of $\mathcal{M}$ by  $Z$ } is the pair
   $\mathcal{M}'=(S\setminus Z,r')$, where  for $A\in [S-Z]^{<\omega }$ we let 
   $$r'(A)=\min\{ r(A\cup Z_{0})-r(Z_{0}):Z_{0}\in [Z]^{<\omega }\}.$$ 
   Since all values of this function are non-negative integers, $r'$ is well-defined (we still need to show that this construction makes a matroid). 

   For $A\in [S-Z]^{<\omega }$, we say that a $Z_{0}\in [Z]^{<\omega }$ {\em fits} $A$ if 
   $$r'(A)=r(A\cup Z_{0})-r(Z_{0}).$$
 \end{definition}
 
 \begin{lemma}\label{MAT_lm:lemma10}
   Let $\mathcal{M}=(S,r)$ be a finitary matroid, $Z\subseteq S$, and $r'$ is the function defined by the contraction.
   
   a) If $A\in [S-Z]^{<\omega }$ and $Z_{0},Z_{1}\in [Z]^{<\omega }$ with $Z_{0}\subseteq Z_{1}$ and $Z_{0}$ fits $A$, then $Z_{1}$ also fits $A$.
   
   b) If $A_{1},...,A_{n}\in [S-Z]^{<\omega }$, then there is a $Z_{0}\in [Z]^{<\omega }$, that fits all of them.
 \end{lemma}
  
 \begin{proof}
  a) In one hand, by the definition of $r'$, we clearly have $r(A\cup Z_{1})-r(Z_{1})\ge r'(A)$. On the other hand, we can write the submodular inequality for $A\cup Z_{0}$ and $Z_{1}$ to get that $r(A\cup Z_{0})+r(Z_{1})\ge r(A\cup Z_{1})+r(Z_{0})$, that can be transformed as $r(A\cup Z_{1})-r(Z_{1})\le r(A\cup Z_{0})-r(Z_{0})=r'(A)$, so the two sides are equal, $Z_{1}$ fits $A$.
  
  b) For each $1\le i\le n$ choose $Z_{i}\in [Z]^{<\omega }$, that fits $A_{i}$. Then by a) $\bigcup_{i=1}^{n}{Z_{i}}$ fits all.
 \end{proof}

\begin{lemma}\label{MAT_lm:lemma11}
   Let $\mathcal{M}=(S,r)$ be a finitary matroid, $Z\subseteq S$, and $r'$ is the function defined by the contraction, then $r'(A)\le r(A)$ for all $A\in [S-Z]^{<\omega }$, and $\mathcal{M'}=(S-Z,r')$ is a finitary matroid.  
 \end{lemma} 
 
\begin{proof}
  By writing $Z_{0}=\emptyset $ into the definition of $r'$, we get that for $A\in [S-Z]^{<\omega }$ $r'(A)\le r(A\cup \emptyset )-r(\emptyset )=r(A)-0=r(A)$. To see that $\mathcal M'$ is a matroid, we need to show properties 1-4 from Definition  \ref{MAT_df:finitary-matroid} for $r'$.
  
  1. Clear by $0\le r'(\emptyset )\le r(\emptyset )=0$
  
  2. Let $A,B\in [S-Z]^{<\omega }$ with $A\subseteq B$ and choose $Z_{0}\in [Z]^{<\omega }$ that fits both. Then since $A\cup Z_{0}\subseteq B\cup Z_{0}$, we have $r(A\cup Z_{0})\le r(B\cup Z_{0})$, so $r'(A)=r(A\cup Z_{0})-r(Z_{0})\le r(B\cup Z_{0})-r(Z_{0})=r'(B)$.
  
  3. Comes from $r'(A)\le r(A)\le |A|$
  
  4. Let $A,B\in [S-Z]^{<\omega }$ and choose $Z_{0}\in [Z]^{<\omega }$, that fits all of $A,B,A\cap B,A\cup B$. Then since $(A\cup Z_{0})\cap (B\cup Z_{0})=(A\cap B)\cup Z_{0}$, and $(A\cup Z_{0})\cup (B\cup Z_{0})=(A\cup B)\cup Z_{0}$, we have $r'(A)+r'(B)=r(A\cup Z_{0})+r(B\cup Z_{0})-2r(Z_{0})\ge r((A\cap B)\cup Z_{0})+r((A\cup B)\cup Z_{0})-2r(Z_{0})=r'(A\cap B)+r'(A\cup B).$
\end{proof} 

\begin{lemma}\label{MAT_lm:lemma12}
   Let $\mathcal{M}=(S,r)$ be a finitary matroid, $Z\subseteq S$, and $\mathcal{M'}=(S-Z,r')$ be the contraction matroid. Then for $X\subseteq S-Z$, we have $X\in \mathcal{I}(\mathcal{M'})$ if and only if for all $Y\in \mathcal{I}(\mathcal{M})\cap \mathcal{P}(Z)$, $X\cup Y\in \mathcal{I}(\mathcal{M})$.
 \end{lemma} 
  
\begin{proof}
  First suppose $X\in \mathcal{I}(\mathcal{M})'$ and let $Y\in \mathcal{I}(\mathcal{M})$ with $Y\subseteq Z$. We need to show that $X\cup Y$ is independent. Let $X_{0}\cup Y_{0}\in [X\cup Y]^{<\omega }$, where $X_{0}\subseteq X$, $Y_{0}\subseteq Y$. Then since $Y$ is independent, we have $r(Y_{0})=|Y_{0}|$ and since $X$ is independent in $\mathcal{M}'$, we have $r(X_{0}\cup Y_{0})-|Y_{0}|=r(X_{0}\cup Y_{0})-r(Y_{0})\ge r'(X_{0})=|X_{0}|$, so $r(X_{0}\cup Y_{0})\ge |X_{0}\cup Y_{0}|$, so they are equal. Since all finite subsets of $X\cup Y$ is in this form, $X\cup Y$ is independent.
  
  Now suppose that $X$ has this property and let $X_{0}\in [X]^{<\omega }$ be arbitrary. Choose a $Z_{0}\in [Z]^{<\omega }$ that fits $X_{0}$, and let $Y\subseteq Z_{0}$ be maximal independent subset of $Z_{0}$, so by lemma \ref{MAT_lm:lemma1}, we have $|Y|=r(Z_{0})$. Since $Y$ is independent, by our assumption $X\cup Y$ is also independent, so $r(X_{0}\cup Y)=|X_{0}\cup Y|=|X_{0}|+|Y|=|X_{0}|+r(Z_{0})$, but then $r'(X_{0})=r(X_{0}\cup Z_{0})-r(Z_{0})\ge r(X_{0}\cup Y)-r(Z_{0})=|X_{0}|+r(Z_{0})-r(Z_{0})=|X_{0}|$. Since $X_{0}$ was arbitrary, we have $X\in \mathcal{I}(\mathcal{M'})$.
\end{proof}
  
\begin{lemma}\label{MAT_lm:lemma13}
   Let $\mathcal{M}=(S,r)$ be a finitary matroid, $Z\subseteq S$. Then the contraction matroid  $\mathcal{M'}=(S-Z,r')$ is loop-free if and only if $Z\subseteq S$ is closed ($\mathcal{M}$ does not need to be loop-free). 
\end{lemma}

\begin{proof}
  First suppose $Z\subseteq S$ is closed. Then for all $x\in S-Z$ and $Z_{0}\in Z^{<\omega }$, we have $r(Z_{0}+x)=r(Z_{0})+1$, so $r(\{ x\}\cup Z_{0})-r(Z_{0})=1$. Then clearly, by the definition, $r'(\{x\})=1$, so $\mathcal{M}'$ is loop-free.
  
  Now suppose that $\mathcal{M}'$ is loop-free. Let $Z_{0}\in [Z]^{<\omega }$ and $x\in S-Z$. Then $r(Z_{0}+x)-r(Z_{0})=r(\{x\} \cup Z_{0})-r(Z_{0})\ge r'(\{x\})=1$, so $r(Z_{0}+x)>r(Z_{0})$. Since $Z_{0}$ and $x$ were arbitrary, we have that $Z$ is closed.
\end{proof}

Now we define the proper colorings of finitary matroids, the chromatic and the list chromatic numbers.

\begin{definition}
   Let $\mathcal{M}=(S,r)$ be a finitary matroid $\mathcal{K}$ be any color set. A 
   {\em proper coloring} of $\mathcal{M}$ by $\mathcal{K}$ is a function $\Phi :S\rightarrow \mathcal{K}$, such that for all $i\in \mathcal{K}$, we have $\Phi ^{-1}(i)\in \mathcal{I}(\mathcal{M})$.
\end{definition}

It is clear that $\Phi $ is a proper coloring if and only if there is no $C\in \mathcal{C}(M)$, such that $C\subseteq \Phi ^{-1}(i)$ for some $i\in \mathcal{K}$.




For loop-free finitary matroids, we can define the chromatic number the following way.

\begin{definition}
  A loop-free finitary matroid $\mathcal{M}=(S,r)$ is {\em $\kappa $-colorable} for a given (finite or infinite) cardinal $\kappa $ if there is a proper coloring $\Phi :S\rightarrow \kappa $ of $\mathcal{M}$. The \emph{chromatic number of $\mathcal{M}$}, denoted by $Chr(\mathcal{M})$, is the smallest cardinal $\kappa$, such that $\mathcal{M}$ is $\kappa $-colorable.
\end{definition}

We can easily see, that for all loop-free matroids $\mathcal{M}=(S,r)$, $Chr(\mathcal{M})$ exists and $Chr(\mathcal{M})\le |S|$.

\begin{definition}
  For a loop-free, finitary matroid $\mathcal{M}=(S,r)$ and a cardinal $\kappa $ a 
  {\em $\kappa $-listing} is a function $L$ from $S$, such that for all $x\in S$, we have $|L(x)|\ge \kappa $. An {\em $L$-coloring} of $\mathcal{M}$ is a function $\Phi :S\rightarrow \bigcup_{x\in S}{L(x)}$, such that $\Phi $ is a proper coloring of $\mathcal{M}$ and for all $x\in S$, we have $\Phi (x)\in L(x)$. $\mathcal{M}$ is {\em $\kappa $ list-colorable}, if it is $L$-colorable for all $\kappa $ listing $L$. The {\em list chromatic number} denoted by $List(\mathcal{M})$ is the smallest cardinal $\kappa $, such that $\mathcal{M}$ is $\kappa $ list-colorable.
\end{definition}

If $\mathcal{M}$ is $\kappa $ list-colorable, then it is also $\kappa $-colorable, as we can write $L(x)=\kappa $ for all $x\in S$, and then the $L$-colorings are exactly the $\kappa $ colorings. By this, if $List(\mathcal{M})$ exists, then $Chr(\mathcal{M})\le List(\mathcal{M})$.

\begin{lemma}\label{MAT_lm:lemma15}
  For all loop-free finitary matroids $\mathcal{M}=(S,r)$ and $List(\mathcal{M})$ exists and $List(\mathcal{M})\le |S|$.
\end{lemma}

\begin{proof}
  Let $\prec $ be a well-ordering of $S$ in order type $|S|$ and let $L$ be any $|S|$-listing. Then by transfinite recursion, for all $x\in S$, since $|L(x)|>|\{\Phi (y):y\prec x\}|$, we can choose $\Phi (x)\in L(x)$, such that $\Phi (x)\neq \Phi (y)$, for $y\prec x$. Then all values of $\Phi $ are different, so it is clearly a proper $L$-coloring.
\end{proof}

Then we have that for all loop-free finitary matroids $Chr(\mathcal{M})\le List(\mathcal{M})\le |S|$. Seymour's famous list-coloring theorem \cite[Theorem 2]{Se98} states that for any finite matroid $Chr(\mathcal{M})=List(\mathcal{M})$ holds. In this paper, we generalize this result for finitary matroid. In the next section \ref{MAT_sc:fin}, we consider  the easy case when $Chr(\mathcal{M})$ is finite, and in Section \ref{MAT_sc:inf} we investigate  the case when $Chr(\mathcal{M})$ is an infinite cardinal.

\section{The finite chromatic number case}\label{MAT_sc:fin}

\begin{theorem}\label{MAT_tm:theorem1}
  Let $\mathcal{M}=(S,r)$ be a loop-free finitary matroid and $k\in \omega $. Then the following statements are equivalent:
  
  (1) $Chr(\mathcal{M})\le k$
  
  (2) $List(\mathcal{M})\le k$
\end{theorem}

\begin{proof}
  $(2)\Rightarrow (1)$ is clear.
  
  $(1)\Rightarrow (2)$
  
  Let $\Phi :S\rightarrow k$ be a $k$-coloring of $\mathcal{M}$. Then clearly for all $A\in [S]^{<\omega }$, $\Phi |_{A}$ is a $k$-coloring of $\mathcal{M}_{A}$, so $Chr(\mathcal{M}_{A})\le k$. Since $\mathcal{M}_{A}$ is a finite matroid, by Seymour's list coloring theorem, we have that $List(\mathcal{M}_{A})\le k$ for all $A\in [S]^{<\omega }$.
  
  Let $L$ be an arbitrary $k$-listing, that we want to construct an $L$-coloring of $\mathcal{M}$. We may suppose that $|L(x)|=k$ for all $x\in S$, as if not, we can choose an arbitrary $L'(x)\subseteq L(x)$ for all $x\in S$ with $|L'(x)|=k$, and then the constructed $L'$-coloring is also an $L$-coloring. For all $A\in [S]^{<\omega }$ choose an $L$-coloring $\Phi _{A}$ of $\mathcal{M}_{A}$.
  
  For any $A\in [S]^{<\omega }$, let $T_{A}=\{ X\in [S]^{<\omega }:A\subseteq S\}$, and let us denote $T_{x}=T_{\{x\}}$. Clearly for any $A$, we have $A\in T_{A}$ and for $A_{1},...,A_{n}\in [S]^{<\omega }$, we have $T_{A_{1}}\cap ...\cap T_{A_{n}}=T_{A_{1}\cup ...\cup A_{n}}$. Since $\mathscr{T}=\{T_{A}:A\in [S]^{<\omega }\}\subseteq \mathcal{P}([S]^{<\omega })$ is a set that is closed under finite intersection and $\emptyset \not\in \mathscr{T}$, 
   there is an ultrafilter $\mathscr{U}\subseteq \mathcal{P}([S]^{<\omega })$, such that $\mathscr{T}\subseteq \mathscr{U}$.
  
  Let $x\in S$ be arbitrary. Then for each $A\in T_{x}$, since $x\in A$, the coloring $\Phi _{A}$ is defined on $x$ and $\Phi _{A}(x)\in L(x)$. For each $i\in L(x)$, let $T_{x,i}=\{A\in T_{x}| \Phi _{A}(x)=i\}$. Then the finite union $\bigcup_{i\in L(x)}{T_{x,i}}=T_{x}\in \mathscr{T}\subseteq \mathscr{U}$, so there is a unique $i\in L(x)$, with $T_{x,i}\in \mathscr{U}$. Define the coloring $\Phi _{\mathscr{U}}:S\rightarrow \bigcup_{x\in S}{L(x)}$, be such that for all $x\in S$, $\Phi _{\mathscr{U}}(x)\in L(x)$ is the unique element with $T_{x,\Phi _{\mathscr{U}}(x)}\in \mathscr{U}$.
  
  We need to show that $\Phi _{\mathscr{U}}$ is a proper coloring. Suppose for contradiction, that there is a $C\in \mathcal{C}(\mathcal{M})$ and some $i$, such that $C\subseteq \Phi _{\mathscr{U}}^{-1}(i)$. Then for all $x\in C$, we have $i\in L(x)$ and $T_{x,i}\in \mathscr{U}$. Then the finite intersection $\bigcap_{x\in C}{T_{x,i}}\in \mathscr{U}$, so $\bigcap_{x\in C}{T_{x,i}}\neq \emptyset $. Let $A\in \bigcap_{x\in C}{T_{x,i}}$. Then for all $x\in C$, $A\in T_{x,i}$, so $x\in A$ and $\Phi _{A}(x)=i$. But then $C\subseteq \Phi _{A}^{-1}(i)$, that is a contradiction, since $\Phi _{A}$ is a proper coloring of $\mathcal{M}_{A}$. Thus, $\Phi _{\mathscr{U}}$ is an $L$-coloring. Hence, $\mathcal{M}$ is $k$ list-colorable, so $List(\mathcal{K})\le k$.
  
\end{proof}
  
 \section{The infinite chromatic number case}\label{MAT_sc:inf}
 
 In this final section, we will show that $Chr(\mathcal{M})=List(\mathcal{M})$, when it is an infinite cardinal. For this, first we need a definition. 
 
 \begin{definition}\label{MAT_df:mb}
  If $\mathcal{M}=(S,r)$ is a loop-free finitary matroid, and 
  we say that $(B,\le )$ is a {\em well-ordered base} of $\mathcal{M}$,
  if $B$ is a base of $\mathcal M$ and $\le $ is a well-order on $B$
  
  Given a well-ordered base $(B,\le )$
  we define the function $M_{B}:S\rightarrow B$ in  the following way:
  \begin{displaymath}
  {M_B(x)}=\left\{\begin{array}{ll}
  {x}&\text{\text{if $x\in B$,}}\\
  {\max_{\le}(C(B,x)\cap B)}&\text{if $x\notin B$,}\\
\end{array}\right.
  \end{displaymath}
 where  $C(B,x)$ denotes the fundamental circuit of $x$ in $B$, see Definition \ref{MAT_df:maincircle}.
  \end{definition}
 
 \begin{theorem}\label{MAT_tm:theorem2} 
  Let $\mathcal{M}=(S,r)$ be a loop-free finitary matroid, $\kappa $ be an infinite cardinal. Then the following statements are equivalent:
  
  (1) $Chr(\mathcal{M})\le \kappa $
  
  (2) $List(\mathcal{M})\le \kappa $
  
  (3) There is a well-ordered base $(B,\le )$ of $\mathcal{M}$, such that for all $b\in B$, $$|\{ x\in S\ |\ M_{B}(x)=b \}|\le \kappa .$$
 \end{theorem}
 
 \begin{proof}
   $(2)\Rightarrow (1)$ is clear.
   
   Before proving  $(3)\Rightarrow (2)$ we need some preparation. Start with an easy observation.
 
 \begin{observation}\label{MAT_lm:lemma16}
  If $\mathcal{M}=(S,r)$ is a finitary matroid $Z\subseteq S$ is closed, and $C\in \mathcal{C}(\mathcal{M})$  such that there is an $x\in C$ with $C-x\subseteq Z$, then $C\subseteq Z$.
 \end{observation}
 
 
 \begin{lemma}\label{MAT_lm:lemma17}
  If $\mathcal{M}=(S,r)$ is a loop-free finitary matroid $(B,\le )$ is a well-ordered base and $C\in \mathcal{C}(\mathcal{M})$. Then there are $x,y\in C$, with $x\neq y$ and $M_{B}(x)=M_{B}(y)$.
 \end{lemma}
 
 \begin{proof}
  Let $B^{*}=(C\cap B)\cup (\bigcup _{x\in C-B}{C(B,x)\cap B})\subset B$. This is a finite subset of $B$. Let $b=max_{\le }(B^{*})$. List the elements of $C$ by $C=\{x_{1},...,x_{n}\}$. Clearly $M_{B}(x)=b$ holds for at least one element of $C$, as either $b\in C$ or $b\in C(B,x)\cap B$ for some $x\in C$, and $b$ is even the maximal element of the whole $B^{*}$. We will show that $M_{b}(x)=b$ holds for at least two elements of $C$. Suppose for contradiction that this is not true. By symmetry, we may suppose, that $M_{B}(x_{n})=b$, and for all $1\le j\le n-1$, $M_{B}(x_{j})\neq b$. Then for $1\le j\le n-1$, if $x_{j}\in B$, then $x_{j}\in B^{*}-b$. If $x_{j}\not\in B$, then $C(B,x_{j})\cap B=C(B,x_{j})-x_{j}\subseteq B^{*}-b$, as if $b$ were in $C(B,x_{j})$, it would be maximal. Thus, by observation \ref{MAT_lm:lemma16}, we have that $x_{1},...,x_{n-1}\in \sigma (B^{*}-b)$. Applying observation  \ref{MAT_lm:lemma16} now for $C$, we also get $x_{n}\in \sigma (B^{*}-b)$. If $x_{n}\in B$, then $x_{n}=M_{B}(x_{n})=b$ would mean that $b\in \sigma (B^{*}-b)$. If $x_{n}\not\in B$, then we have $C(B,x_{n})-x_{n}-b\subseteq B^{*}-b$ and $x_{n}\in \sigma (B^{*}-b)$, so $C(B,x_{n})-b\subseteq \sigma (B^{*}-b)$, using observation  \ref{MAT_lm:lemma16} again, we also get that $b\in \sigma (B^{*}-b)$. But since $B^{*}$ is independent, for all $B_{0}\in [B^{*}-b]^{<\omega }$, we have $r(B_{0}+b)=|B_{0}+b|=r(B_{0})+1>r(B_{0})$ in contradiction with lemma  \ref{MAT_lm:lemma8}. Thus, $M_{b}(x)=b$ holds for at least two $x\in C$.
 \end{proof}
   $(3)\Rightarrow (2)$
 
 Suppose, that (3) holds, and let $(B,\le )$ be a well-ordered basis of $\mathcal{M}$ such that for all $b\in B$, $|S_{b}|=|\{x\in S|M_{B}(x)=b\}|\le \kappa $. We need to show  that $\mathcal{M}$ is $\kappa $ list-colorable. Let $L$ be an arbitrary $\kappa $-listing. First, we will show  that we can choose a $\Phi _{b}$ function on $S_{b}$, such that for $x\in S_{b}$, we have $\Phi _{b}(x)\in L(x)$ and $\Phi _{b}$ is one-to-one. Let $\prec _{b}$ be a well-ordering of $S_{b}$ in order type $\le \kappa $ and let us define $\Phi _{b}$ by transfinite recursion. For $x\in S_{b}$, since $|L(x)|\ge \kappa $ and $|\{\Phi_{b}(y),y\prec _{b}x\}|<\kappa $, we can choose $\Phi _{b}(x)\in L(x)$, such that $\Phi _{b}(x)\neq \Phi _{b}(y)$ for $y\prec _{b}x$. The $\Phi _{b}$ defined this way is clearly one-to-one. Let $\Phi =\bigcup_{b\in B}{\Phi _{b}}$. Then clearly for all $x\in S$, we have $\Phi (x)=\Phi_{M_{B}(x)}(x)\in L(x)$. We need to show  that this is a proper coloring. Suppose for contradiction, that there is some $C\in \mathcal{C}(\mathcal{M})$ and an $i\in \bigcup_{x\in S}{L(x)}$, such that $C\subseteq \Phi^{-1}(i)$. Then by lemma  \ref{MAT_lm:lemma17}, there are $x,y\in C$, $x\neq y$ and $M_{B}(x)=M_{B}(y)=b\in B$. Then $x,y\in S_{b}$, so $\Phi (x)=\Phi _{b}(x)\neq \Phi _{b}(y)=\Phi (y)$, since $\Phi _{b}$ is one-to-one. This is in contradiction with $\Phi (x)=i=\Phi (y)$, so $\Phi $ is a proper coloring.
 
Before proving $(1)\Rightarrow (3)$ we need to prove some lemmas.

 \begin{lemma}\label{MAT_lm:lemma18}
   Let $\mathcal{M}=(S,r)$ be a finitary matroid, $A\in [S]^{<\omega }$, with $|A|=n$, and $x_{1},...,x_{n+1}\in S-A$ be such that $r(A+x_{i})=r(A)$ for all $1\le i\le n+1$. Then the set $\{x_{1},...,x_{n+1}\}$ is not independent.
 \end{lemma}
 
 \begin{proof}
   By lemma  \ref{MAT_lm:lemma5}, we have that $r(A\cup \{x_{1},...,x_{n+1}\})=r(A)\le |A|=n$, so by monotonity $r(\{x_{1},...,x_{n+1}\})\le r(A\cup \{x_{1},...,x_{n+1}\})<n+1=|\{x_{1},...,x_{n+1}\}|$, that it is not independent.
 \end{proof}
  
  \begin{lemma}\label{MAT_lm:lemma19}
   Let $\mathcal{M}=(S,r)$ be a loop-free finitary matroid, $\kappa $ be an infinite cardinal, with $Chr(\mathcal{M})\le \kappa $. Then for all $[A]\in [S]^{<\omega }$, we have 
   $$|\{x\in S\setminus A\ |\ r(A+x)=r(A) \}|\le \kappa .$$
  \end{lemma}
  
  \begin{proof}
   Let $\Phi :S\rightarrow \kappa $ be a $\kappa $-coloring of $\mathcal{M}$ and for all $\alpha <\kappa $ let $A_{\alpha }=\{ x\in S-A|r(A+x)=r(A),\Phi (x)=\alpha \}$. Let $n=|A|$. First we will show that $|A_{\alpha }|\le n$ for all $\alpha $. Suppose for contradiction, that for some $\alpha $, we have $|A_{\alpha }|\ge n+1$, and let $x_{1},...,x_{n+1}\in A_{\alpha}$. Then by lemma  \ref{MAT_lm:lemma18}, we have that $\{x_{1},...,x_{n+1}\}\not\in \mathcal{I}(\mathcal{M})$ and $\{x_{1},...,x_{n+1}\}\subseteq A_{\alpha}\subseteq \Phi^{-1}(\alpha )$ in contradiction with $\Phi $ is a proper coloring. Then we have $|\{x\in S-A|r(A+x)=r(A) \}|=|\bigcup_{\alpha <\kappa }{A_{\alpha }}|\le \kappa \cdot n=\kappa $.
  \end{proof}

  Thus, by Lemma \ref{MAT_lm:lemma19}, for each 
   for loop-free finitary matroid $\mathcal M$ with $Chr(\mathcal{M})\le \kappa $
 we can fix  a {\em bookkeeping function $h$}, i.e. 
   a function $h:[S]^{<\omega }\times \kappa \rightarrow S$ such that for all $A\in [S]^{<\omega }$
   \begin{displaymath}
    \{x\in S: r(A+x)=r(A)\}\subset \{h(A,{\alpha}):{\alpha}<{\kappa}\}.
   \end{displaymath}
  
  The last thing, we need for this proof is some model theoretic approach, using elementary submodels. For this, let $H(\theta )=\{ x:|TC(x)|<\theta \}$ for some infinite cardinal $\theta $. Here $TC$ denotes the transitive closure, so $TC(x)=\bigcup_{n\in \omega }{U_{n}(x)}$, where $U_{0}(x)=x$ and $U_{n+1}(x)=\cup U_{n}(x)$, for any set $x$. If $\theta $ is a regular cardinal, then $H(\theta )$ is a set and all ZFC axioms except Power Set Axiom holds in the model $H(\theta )$. However, if $x\in H(\theta )$ and $2^{|x|}<\theta $, then we also have $\mathcal{P}(X)\in H(\theta )$. Since all sets belong to some $H(\theta )$, in practice we may say that $\theta $ is sufficiently large. That means, that all sets defined in the proof are in $H(\theta )$.
   An $M$ is an elementary submodel of $H(\theta )$ if for set-theoretic formulae $\phi $ and $x_{1},...,x_{n}\in M$ (where $n$ is the number of free variables of $\phi $) we have $M\vDash \phi (x_{1},...,x_{n})\Leftrightarrow H(\theta ){\vDash}\phi (x_{1},...,x_{n})$. By Löwenheim-Skolem theorem, for all $R\subseteq H(\theta )$, there is an elementary submodel $M$, such that $R\subseteq M$ and $|M|=max(|R|,\omega )$.
  
  $(1)\Rightarrow (3)$
  
  Let $\kappa $ be fixed. We define the statement $Q(\lambda )$ for each infinite cardinal $\lambda $ in  the following way: 
  \begin{enumerate}[($Q({\lambda})$)]
  \item If $\mathcal{M}=(S,r)$ is a loop-free matroid, with $Chr(\mathcal{M})\le \kappa $ and $|S|=\lambda $, then there is a well-ordered base $(B,\le )$ of $\mathcal{M}$, such that for all $b\in B$, we have $|\{x\in S\ |\ M_{B}(x)=b\}|\le \kappa $. 
  \end{enumerate}
  
  The statement $(1)\Rightarrow (3)$ would mean, that $Q(\lambda )$ holds for all cardinals $\lambda $. We prove it by induction on $\lambda $.
  
  If $\lambda \le \kappa $, then $Q(\lambda )$ clearly holds, as any well-ordered base $(B,\le )$ would fit. Now suppose that $\lambda >\kappa $ and $Q(\mu )$ holds for all $\mu <\lambda $. We need to prove that $Q(\lambda )$ also holds.
  
  Let $\mathcal{M}=(S,r)$ be a finitary matroid with $|S|=\lambda $ and $Chr(\mathcal{M})\le \kappa $. Fix a proper coloring $\Phi :S\rightarrow \kappa $ and a bookkeeping function $h:[S]^{<\omega }\times \kappa \rightarrow S$. Let $(S_{\alpha})_{\alpha <cf(\lambda)}$ be such that $S_{\alpha }\subseteq S$, $|S_{\alpha} |<\lambda$ for all $\alpha <cf(\lambda )$ and $\bigcup _{\alpha <cf(\lambda )}{S_{\alpha }}=S$. Let $\theta $ be a sufficiently large regular cardinal. We construct an increasing sequence of elementary submodels $M_{\alpha }\subseteq H(\theta )$ for $\alpha <cf(\lambda )$. Let $M_{0}\subseteq H(\theta )$ be an elementary submodel, such that $\kappa \cup \{\mathcal{M},\Phi ,h\}\subseteq M_{0}$ and $|M_{0}|=\kappa $. For $\alpha <cf(\lambda )$, let $M_{\alpha +1}\subseteq H(\theta )$ be an elementary submodel with $M_{\alpha }\cup S_{\alpha }\subseteq M_{\alpha +1}$, and $|M_{\alpha +1}|=|M_{\alpha }\cup S_{\alpha }|$. If $\alpha <cf(\lambda )$ is a limit ordinal, then let $M_{\alpha }=\bigcup_{\beta <\alpha }{M_{\beta}}$. 
  We also have that for all $\alpha <cf(\lambda )$, $|M_{\alpha }|<\lambda $. For $M_{0}$, we have $|M_{0}|=\kappa <\lambda $, for successor ordinals, since $|M_{\alpha }|<\lambda $ and $|S_{\alpha }|<\lambda $, we have by definition, that $|M_{\alpha +1}|<\lambda $. For limit ordinals $M_{\alpha }$ is a $<cf(\lambda )$ union of $<\lambda $ sets, so $|M_{\alpha }|<\lambda $ also holds. For all $\alpha <cf(\lambda )$, let $Z_{\alpha }=M_{\alpha }\cap S$.
  
  \begin{lemma}\label{MAT_lm:lemma20}
   For all $\alpha <cf(\lambda )$, the set $Z_{\alpha }\subseteq S$ is closed.
  \end{lemma}
  
  \begin{proof}
   Suppose for contradiction, that $Z_{\alpha }$ is not closed and let $A\in [Z_{\alpha}]^{<\omega }$, and $x\in S-Z_{\alpha }$ be such that $r(A+x)=r(A)$. Then since $A\subseteq Z_{\alpha }\subseteq M_{\alpha}$ is a finite subset of an elementary submodel, we have $A\in M_{\alpha }$. Since $h$ is a bookkeeping function, there is some $\gamma <\kappa $, such that $h(A,\gamma )=x$. Clearly, we also have $h,\gamma \in M_{0}\subseteq M_{\alpha }$, so $x=h(A,\gamma )\in M_{\alpha }$. Thus, $x\in Z_{\alpha }$, that is a contradiction.
  \end{proof}
  
  Now since for all $\alpha <cf(\lambda )$, we have $S_{\alpha }\subseteq Z_{\alpha +1}$ and $S=\bigcup _{\alpha <cf(\lambda )}{S_{\alpha }}$, for all $x\in S$ there is some $\alpha <cf(\lambda )$, such that $x\in Z_{\alpha }$. Let us define the rank function $\rho :S\rightarrow cf(\lambda )$, such that $\rho (x)=min\{\alpha :x\in Z_{\alpha}\}$. Moreover, $\rho (x)$ must be a successor ordinal, as for limit ordinals $Z_{\alpha }$ is just the union of former ones.
   
  We construct an increasing chain $(B_{\alpha },\le _{\alpha })_{\alpha <cf(\lambda )}$ of well-ordered independent sets in $\mathcal{M}$, such that 
  \begin{enumerate}[(i)]
  \item $B_{\alpha }\in \mathcal{B}(\mathcal{M}_{Z_{\alpha }})$,
  \item  
  $(B_{\alpha},\le_{\alpha})$ is an initial segment of $(B_{\beta},\le_{\beta})$ for 
  ${\alpha}<{\beta}<cf({\lambda})$,
  \item 
   $|\{x\in Z_{\alpha }\ |\ M_{B_{\alpha }}(x)=b\}|\le \kappa $ for each ${\alpha}<cf({\lambda})$ and $b\in B_{\alpha }$. 
  \end{enumerate}
  We construct those $B_{\alpha}$s  by transfinite recursion. First, since $Q(\kappa )$ clearly holds and $|Z_{0}|\le |M_{0}|=\kappa $, we can construct $(B_{0},\le _{0})$.
  
  Now suppose, that $\alpha $ is a limit ordinal and for $\beta <\alpha $, $(B_{\beta },\le _{\beta })$ is already constructed. Then clearly $Z_{\alpha }=M_{\alpha }\cap S=(\bigcup_{\beta <\alpha }{M_{\beta}})\cap S=\bigcup_{\beta <\alpha }{(M_{\beta}\cap S)}=\bigcup_{\beta <\alpha }{Z_{\beta }}$. Let $B_{\alpha }=\bigcup _{\beta <\alpha }{B_{\beta }}$
  and $\le_{\alpha }=\bigcup _{\beta <\alpha }{\le_{\beta }}$
  
  We need to show that $B_{\alpha }\in \mathcal{B}(\mathcal{M}_{Z_{\alpha }})$. First we show that $B_{\alpha }$ is independent. Let $B'\in [B_{\alpha }]^{<\omega }$ and let $B'=\{x_{1},...,x_{n}\}$. Since for all $1\le i\le n$, we have $x_{i}\in B_{\alpha }=\bigcup _{\beta <\alpha }{B_{\beta }}$, there is some $\beta _{i}<\alpha $, such that $x_{i}\in B_{\beta _{i}}$. Let $\beta =max_{1\le i\le n}(\beta _{i})$. Then for all $i$, $x_{i}\in B_{\beta _{i}}\subseteq B_{\beta }$, so $B'\subseteq B_{\beta }$. Since $B_{\beta }$ is independent, we have $r(B')=|B'|=n$. Since it was an arbitrary finite subset, we have $B_{\alpha }\in \mathcal{I}(\mathcal{M}_{Z_{\alpha }})$. We also show  that this is maximal. Let $x\in Z_{\alpha }-B_{\alpha }$ be arbitrary. Then for $\beta =\rho (x)<\alpha $, since $x\in Z _{\beta }$, $B_{\beta }+x$ is not independent, thus $B_{\alpha }+x$ neither. Hence, $B_{\alpha }$ is maximal independent in $Z_{\alpha }$, so it is a base.

  By (ii), $\le_{\alpha}$ is a well-ordering of $B_{\alpha}$,
  and $(B_{\beta},\le_{\beta})$ is an initial segment of $(B_{\alpha},\le_{\alpha})$.


  For (iii), let $b\in B_{\alpha }$ be arbitrary, and let 
  $\beta =\min\{{\gamma}<{\alpha}:b\in Z_{\gamma}\}<\alpha 
  $. Now, for any $x\in Z_{\alpha }$ with $M_{B_{\alpha }}(x)=b\in B_{\beta }$, we have either $x=b$ or $x\not \in B_{\alpha }$ and $C(B_{\alpha },x)-x=C(B_{\alpha },x)\cap B_{\alpha }\subseteq B_{\beta}\subseteq Z_{\beta }$, as $B_{\beta }$ is an initial segment. Then by lemma  \ref{MAT_lm:lemma20}, $Z_{\alpha }$ is closed, applying observation  \ref{MAT_lm:lemma16}, we get that $x\in B_{\beta }$. Here we also have that $M_{B_{\alpha }}(x)=M_{B_{\beta }}(x)$. If $x\in B_{\beta }\subseteq B_{\alpha }$, this is clear, otherwise $x\not\in B_{\beta }$, so $C(B_{\beta },x)\subseteq B_{\beta }+x\subseteq B_{\alpha }+x$, that is a circuit, so by lemma  \ref{MAT_lm:lemma3}, we have that $C(B_{\alpha },x)=C(B_{\beta },x)$, thus $M_{B_{\alpha }}(x)=M_{B_{\beta }}(x)$. Hence, $|\{x\in Z_{\alpha }|M_{B_{\alpha }}(x)=b\}|=|\{x\in Z_{\beta }|M_{B_{\beta }}(x)=b\}|\le \kappa $, so we are done.
  
  \medskip

  Now the successor case: let $\alpha <cf(\lambda )$ and assume that 
  $(B_{\alpha},\le_{\alpha})$ is constructed. In the matroid $\mathcal{M}_{Z_{\alpha +1}}$, the set $Z_{\alpha }$ is closed by lemma  \ref{MAT_lm:lemma20}. Then, by lemma  \ref{MAT_lm:lemma13}, the contracted matroid $\mathcal{M'}=(Z_{\alpha +1}-Z_{\alpha },r')$ is loop-free. Next we will see that $Chr(\mathcal{M}')\le \kappa $ also holds.
  
  \begin{lemma}\label{MAT_lm:lemma21}
   For the contracted matroid $\mathcal{M'}=(Z_{\alpha +1}-Z_{\alpha },r')$, the restriction $\Phi |_{Z_{\alpha +1}-Z_{\alpha }}$ is a proper coloring of $\mathcal{M}'$, and thus $Chr(\mathcal{M}')\le \kappa $. 
  \end{lemma}
  
  \begin{proof}
   Suppose for contradiction, that $\Phi $ is not a proper coloring. Then there is an $X\subseteq Z_{\alpha +1}-Z_{\alpha }$ and a $\gamma <\kappa $, such that $X\subseteq \Phi ^{-1}(\gamma )$ and $X\not\in \mathcal{I}(\mathcal{M'})$. 
   Then by lemma  \ref{MAT_lm:lemma12}, there is a $Y\subseteq Z_{\alpha }$, such that $Y\in \mathcal{I}(\mathcal{M})$, but $X\cup Y\not\in \mathcal{I}(\mathcal{M})$. Let $C\in \mathcal{C}(\mathcal{M})$, be such that $C\subseteq X\cup Y$. 
   Then since $Y$ is independent, we must have $C\not\subseteq Y$, and since $C\cap Z_{\alpha }\subseteq Y$, we have $C\not\subseteq Z_{\alpha }$. Moreover, for all $x\in C-Z_{\alpha }$, we have $x\in X\subseteq \Phi ^{-1}(\gamma )$, so $\Phi (x)=\gamma $.
   
   Let 
   \begin{multline*}
    k=min \{ |A|:A\in [Z_{\alpha }]^{<\omega }\land 
    \exists C'\in \mathcal{C}(\mathcal{M}_{Z_{\alpha +1}})
    \\ C'\not\subseteq Z_{\alpha }\land \forall x\in C'\setminus A\ ,\Phi (x)=\gamma \}.
   \end{multline*}
   
   As we can place $A_{0}=C-\Phi ^{-1}(\gamma )\subseteq Z_{\alpha }$ into this definition, so $k$ is well-defined. We also must have $k\ge 1$, as for $k=0$, we would have $C'\subseteq \Phi^{-1}(\gamma )$ in contradiction with $\Phi $ is a proper coloring of $\mathcal M$. Let $A\in [Z_{\alpha }]^{<\omega }$, be a set with $|A|=k$, and $C_{1}\in \mathcal{C}(\mathcal{M}_{Z_{\alpha +1 }})$ be a circuit with $A\subseteq C_{1}$,  such that $C_{1}\not \subseteq Z_{\alpha }$ and $C_{1}-A\subseteq \Phi ^{-1}(\gamma )$. Let $l=|C_{1}|-k=|C_{1}-A|$. 
   Then we have that $$H(\theta )\vDash\exists x_{1}...\exists x_{l}, \Phi (x_{1})=\gamma \wedge ...\wedge \Phi (x_{l})=\gamma \wedge A\cup \{x_{1},...,x_{l}\}\in \mathcal{C}(\mathcal{M}).$$ 
   Since $A\subseteq Z_{\alpha }\subseteq M_{\alpha }$ is a finite subset, we have $A\in M_{\alpha }$. As $\Phi ,\gamma , \mathcal{M}\in M_{0}\subseteq M_{\alpha }$ and  $M_{\alpha }$ is an elementary submodel of $H(\theta )$, we have 
   $$M_{\alpha }\vDash \exists x_{1}...\exists x_{l},\Phi (x_{1})=\gamma \wedge ...\wedge \Phi (x_{l})=\gamma \wedge A\cup \{x_{1},...,x_{l}\}\in \mathcal{C}(\mathcal{M}).$$ 
   Then there is a $C_{2}\in \mathcal{C}(\mathcal{M})$, with $A\subseteq C_{2}\subseteq Z_{\alpha }$ and for all $x\in C_{2}-A$, $\Phi (x)=\gamma $. Let $e\in A$, then clearly $e\in C_{1}\cap C_{2}$. Also choose an $e_{1}\in C_{1}-Z_{\alpha }\subseteq C_{1}-C_{2}$. Then we can use lemma  \ref{MAT_lm:lemma2} b), so there is a $C_{3}\in \mathcal{C}(\mathcal{M})$, with $C_{3}\subseteq (C_{1}\cup C_{2})-e$ and $e_{1}\in C_{3}$. 
   Then clearly $C_{3}\subseteq C_{1}\cup C_{2}\subseteq Z_{\alpha +1}$, and $C_{3}\not\subseteq Z_{\alpha }$, as $e_{1}\in C_{3}$. Let $A_{1}=A\cap C_{3}$. Since $A_{1}\subseteq A-e$, $|A_{1}|\le k-1$. 
   Moreover, for all $x\in C_{3}-A_{1}$, we either have $x\in C_{1}-A$, or $x\in C_{2}-A$, thus $\Phi (x)=\gamma $. 
   Then the set $A_{1}$ is also good with the circuit $C_{3}$, that is a contradiction with the minimality of $k$. Hence, $\Phi |_{Z_{\alpha +1}-Z_{\alpha }}$ is a proper coloring of $\mathcal{M}'$.
  \end{proof}
  
  Now let $\mu =|Z_{\alpha +1}-Z_{\alpha }|\le |Z_{\alpha +1}|\le |M_{\alpha +1}|<\lambda $. By $Q(\mu )$, $\mathcal{M}'$  has a well ordered base $(B_{\alpha }',\le _{\alpha }')$, such that  $|\{x\in Z_{\alpha +1}\setminus Z_{\alpha }:M_{B_{\alpha }'}(x)=b\}|\le \kappa $
  for each $b\in B_{\alpha }'$. 
  
  Let $B_{\alpha +1}=B_{\alpha }\cup B_{\alpha }^{'}$. First we need to show  that this is a base. 
  For this, we prove the properties of lemma  \ref{MAT_lm:lemma9}. Since $B_{\alpha}$ is independent, and $B_{\alpha }'$ is independent in the contracted matroid, by lemma  \ref{MAT_lm:lemma12}, we have $B_{\alpha +1}$ is also independent. 
  Now we need to show that $\sigma (B_{\alpha +1})=Z_{\alpha +1}$. Since by lemma  \ref{MAT_lm:lemma20} $Z_{\alpha +1}$ is closed, we have $\sigma (B_{\alpha +1})\subseteq Z_{\alpha +1}$. Let $x\in Z_{\alpha +1}$. If $x\in Z_{\alpha }$, then using lemma  \ref{MAT_lm:lemma9}, we get $x\in \sigma (B_{\alpha })\subseteq \sigma (B_{\alpha +1})$. 
  Suppose $x\in Z_{\alpha +1}-Z_{\alpha }$. If $x\in B_{\alpha }'$, then we are done, so suppose $x\not\in B_{\alpha }'$. Then $B_{\alpha }'+x$ is not independent in $\mathcal{M}'$, so by lemma  \ref{MAT_lm:lemma12}, there is an independent $Y\subseteq Z_{\alpha }$, such that $(Y\cup B_{\alpha }')+x$ is not independent. Let $C\in \mathcal{C}(\mathcal{M})$, with $C\subseteq (Y\cup B_{\alpha }')+x$. 
  Since $B_{\alpha }'\in \mathcal{I}(\mathcal{M})'$, by lemma  \ref{MAT_lm:lemma12}, we have $Y\cup B_{\alpha }'\in \mathcal{I}(\mathcal{M})$, we must have $x\in C$. 
  In one hand, we have $Y\subseteq Z_{\alpha } =\sigma (B_{\alpha })\subseteq \sigma (B_{\alpha +1})$, on the other hand $B_{\alpha }'\subseteq B_{\alpha +1}\subseteq \sigma (B_{\alpha +1})$, thus $C-x\subseteq Y\cup B_{\alpha }'\subseteq \sigma (B_{\alpha +1})$, so by observation  \ref{MAT_lm:lemma16}, we have $x\in \sigma (B_{\alpha +1})$. Since $x\in Z_{\alpha +1}$ was arbitrary, $\sigma (B_{\alpha +1})=Z_{\alpha +1}$ by lemma  \ref{MAT_lm:lemma9} $B_{\alpha +1}$ is a base.
  
  For the well ordering $\le _{\alpha +1}$, we simply put $B_{\alpha }'$ into the top of $B_{\alpha }$, more formally, 
  \begin{displaymath}
  \le_{{\alpha}+1}\ =\ \le_{\alpha}\cup (B_{\alpha}\times B'_{\alpha})\cup \le'_{\alpha}.
  \end{displaymath}
  
  
  Clearly this is a well ordering and $B_{\alpha }$ (and hence all $B_{\beta }$-s for $\beta <\alpha $ )is an initial segment. Now we need to show that for all $b\in B_{\alpha +1}$, we have $|\{x\in Z_{\alpha +1}|M_{B_{\alpha +1}}(x)=b\}|\le \kappa $. First suppose $b\in B_{\alpha }$. We will show  that for any $x\in Z_{\alpha +1}$ with $M_{B_{\alpha +1}}(x)=b$, we have $x\in Z_{\alpha }$ and $M_{B_{\alpha }}(x)=b$. If $x\in B_{\alpha +1}$, then $x=b$, so this is clear. Suppose that $x\not\in B_{\alpha +1}$. Then since $B_{\alpha }$ is an initial segment of $B_{\alpha +1}$, we have that $C(B_{\alpha +1},x)-x=C(B_{\alpha +1},x)\cap B_{\alpha +1}\subseteq B_{\alpha }\subseteq Z_{\alpha }$. By lemma  \ref{MAT_lm:lemma20}, $Z_{\alpha }$ is closed, so by observation  \ref{MAT_lm:lemma16}, we have $x\in Z_{\alpha }$. Moreover, $C(B_{\alpha },x)\subseteq B_{\alpha }+x\subseteq B_{\alpha +1}+x$ is a circuit, so by lemma  \ref{MAT_lm:lemma3}, $C(B_{\alpha +1},x)=C(B_{\alpha },x)$, thus $b=M_{B_{\alpha +1}}(x)=M_{B_{\alpha }}(x)$. Hence, $|\{x\in Z_{\alpha +1}|M_{B_{\alpha +1}}(x)=b\}|=|\{x\in Z_{\alpha }|M_{B_{\alpha }}(x)=b\}|\le \kappa $. Now suppose $b\in B_{\alpha }'$. 
  We will show  that for any $x\in Z_{\alpha +1}$ with $M_{B_{\alpha +1}}(x)=b$, we have $x\in Z_{\alpha +1}-Z_{\alpha }$ and $M_{B_{\alpha }'}(x)=b$. Again if $x\in B_{\alpha +1}$, then $x=b$, so this is clear, so suppose that $x\not\in B_{\alpha +1}$. If $x$ was in $Z_{\alpha }$, we would have a circuit $C(B_{\alpha },x)\subseteq B_{\alpha }+x\subseteq B_{\alpha +1}+x$, so by lemma  \ref{MAT_lm:lemma3}, we would have $C(B_{\alpha +1},x)=C(B_{\alpha },x)$, thus $b=M_{B_{\alpha +1}}(x)=M_{B_{\alpha }}(x)$, that is a contradiction. So we must have $x\in Z_{\alpha +1}-Z_{\alpha }$. Let $C'=C(B_{\alpha +1},x)-Z_{\alpha }$. Then $C'\subseteq B_{\alpha '}+x$. We will show  that $C'$ is the fundamental circuit of $x$ on $B_{\alpha }'$ in the matroid $\mathcal{M}'$. Then by lemma  \ref{MAT_lm:lemma3}, we only need to prove, that it is a circuit in the contraction matroid. First $C'$ is not independent in $\mathcal{M}'$, as if it was independent, then by lemma  \ref{MAT_lm:lemma12}, $C(B_{\alpha +1},x)=C'\cup (C(B_{\alpha +1},x)\cap Z_{\alpha })=C'\cup (C(B_{\alpha +1},x)\cap B_{\alpha })$ would be independent in $\mathcal{M}$, that is not true. Now we need to show  that for all $X\subset C'$ with $X\neq C'$, we have $X\in \mathcal{I}(\mathcal{M}')$. If $x\not\in X$, then $X\subseteq B_{\alpha }'\in \mathcal{I}(\mathcal{M}')$. Suppose $x\in X$, and suppose for contradiction, that $X\not\in \mathcal{I(\mathcal{M}')}$ Then by lemma  \ref{MAT_lm:lemma12}, there is a $Y\subseteq Z_{\alpha }$ independent set, such that $X\cup Y$ is not independent. Let $C\in \mathcal{C}(\mathcal{M})$ be such that $\mathcal{C}\subseteq X\cup Y$. As $X-x\subseteq B_{\alpha }'\in \mathcal{I}(\mathcal{M})'$, by lemma  \ref{MAT_lm:lemma12} we have $(X-x)\cup Y\in \mathcal{I}(\mathcal{M})$, so we must have $x\in C$. Then by lemma  \ref{MAT_lm:lemma9}, we have $Y\subseteq Z_{\alpha }=\sigma (B_{\alpha })\subseteq \sigma (B_{\alpha }\cup (X-x))$, thus $C-x\subseteq Y\cup (X-x)\subseteq \sigma (B_{\alpha }\cup (X-x))$, so by observation  \ref{MAT_lm:lemma16}, we have $x\in \sigma (B_{\alpha }\cup (X-x))$. Then by lemma  \ref{MAT_lm:lemma8}, there is an $A\in [B_{\alpha }\cup (X-x)]^{<\omega }$, such that $r(A+x)=r(A)\le |A|<|A+x|$, thus $A+x$ is not independent. Let $C''\in \mathcal{C}(\mathcal{M})$ be such that $C''\subseteq A+x$. But then $C''\subseteq A+x\subseteq B_{\alpha }\cup X\subseteq B_{\alpha +1}+x$, then by lemma  \ref{MAT_lm:lemma3}, we have $C''=C(B_{\alpha +1},x)$. But $C''-Z_{\alpha }\subseteq X$, witch is a proper subset of $C'$ and  $C(B_{\alpha +1},x)-Z_{\alpha }=C'$, that is a contradiction. Hence, $C'$ is a circuit in $\mathcal{M}'$, so it is the fundamental circuit of $x$ for $B_{\alpha }'$. Then $b=M_{B_{\alpha +1}}(x)=max_{\le _{\alpha }'}(C'-x)=M_{B_{\alpha }'}(x)$. Thus, $|\{x\in Z_{\alpha +1}|M_{B_{\alpha +1}}(x)=b\}|=|\{x\in Z_{\alpha }|M_{B_{\alpha }'}(x)=b\}|\le \kappa $. This way the induction step works for successor cardinals.
  
  Finally, let $B=\bigcup _{\alpha <cf(\lambda )}{B_{\alpha }}$, and   we define the well ordering $\le $ on $B$ by taking $\le=\bigcup _{\alpha <cf(\lambda )}{\le _{\alpha }}$
   Similarly as in the proof for the limit step,   we can see that  
   $|\{x\in S|M_{B}(x)=b\}|\le \kappa $, thus $Q(\lambda )$ for all $b\in B$. 
   
   Thus, by transfinite induction, we have proved, that $Q(\lambda )$ holds for all cardinals, so $(1)\Rightarrow (3)$ holds. 
 \end{proof}
 
 \begin{theorem}\label{MAT_tm:theorem3}
 For any loop-free finitary matroid $\mathcal{M}=(S,r)$, we have $Chr(\mathcal{M})=List(\mathcal{M})$.
 \end{theorem}
 
 \begin{proof} 
 If $Chr(\mathcal{M})=k\in \omega $, by theorem  \ref{MAT_tm:theorem1}, we have $List(\mathcal{M})\le k=Chr(\mathcal{M})$. If $Chr(\mathcal{M})$ is an infinite cardinal, then by theorem \ref{MAT_tm:theorem2} $List(\mathcal{M})\le k=Chr(\mathcal{M})$. As clearly $Chr(\mathcal{M})\le List(\mathcal{M})$, $Chr(\mathcal{M})= List(\mathcal{M})$ holds for all loop-free finitary matroids.
 \end{proof}

\section{Cooperative coloring of finitary matroids}\label{MAT:sc_coop}

In this section we will show the following theorem:

\begin{theorem}\label{TH_CooperativeFinitaryGeneral}
    Let $S$ be a base set, $\mathcal{K}$ be a color set $(\mathcal{M}_{i})_{i\in \mathcal{K}}$ be a collection of loop-free finitary matroids on the set $S$ and $L:S\rightarrow \mathcal{P}(\mathcal{K})$ be a listing. Suppose that there is a (finite or infinite) cardinal $\kappa $, such that $Chr(\mathcal{M}_{i})\le \kappa $ for all $i\in \mathcal{K}$ and $|L(x)|\ge \kappa $ for all $x\in S$. Then there is a function $\Phi :S\rightarrow \mathcal{K}$, such that

    1. $\Phi(x)\in L(x)$ for all $x\in S$

    2. $\Phi ^{-1}(i)\in \mathcal{I}(\mathcal{M}_{i})$ for all $i\in \mathcal{K}$.
\end{theorem}

\begin{remark}\label{RE_SeymourListColoring}
    Theorem \ref{TH_CooperativeFinitaryGeneral} clearly shows the earlier result $Chr(\mathcal{M})=List(\mathcal{M})$ of theorem \ref{MAT_tm:theorem3} as for any listing $L$, we can apply this theorem with $\mathcal{K}=\bigcup_{x\in S}{L(x)}$, $\mathcal{M}_{i}=\mathcal{M}$ for all $i\in \mathcal{K}$ and $\kappa =Chr(\mathcal{M})$.
\end{remark}

\begin{proof} (Of theorem \ref{TH_CooperativeFinitaryGeneral})

First we will show it for $\kappa =k\in \omega $

\begin{lemma}\label{LM_CooperativeFinitaryFinite}
    Let $S$ be a base set, $\mathcal{K}$ be a color set $(\mathcal{M}_{i})_{i\in \mathcal{K}}$ be a collection of loop-free finitary matroids on the set $S$ and $L:S\rightarrow \mathcal{P}(\mathcal{K})$ be a listing. Suppose that there is $k\in \omega $, such that $Chr(\mathcal{M}_{i})\le k$ for all $i\in \mathcal{K}$ and $|L(x)|\ge k $ for all $x\in S$. Then there is a function $\Phi :S\rightarrow \mathcal{K}$, such that

    1. $\Phi(x)\in L(x)$ for all $x\in S$

    2. $\Phi ^{-1}(i)\in \mathcal{I}(\mathcal{M}_{i})$ for all $i\in \mathcal{K}$.
\end{lemma}

\begin{proof}
    First, we may suppose, that $|L(x)|=k$ for all $x\in S$, otherwise we can replace them with $L'(x)\subseteq L(x)$ with $|L'(x)|=k$.

    For any $A\in [S]^{<\omega }$, let $\mathcal{K}_{A}=\bigcup_{x\in A}{L(x)}$. This is clearly a finite set and for all $i\in \mathcal{K}_{A}$, let $\mathcal{M}_{i,A}$ be the matroid $\mathcal{M}_{i}$ restricted to the set $A$. These are clearly finite matroids with $Chr(\mathcal{M}_{i,A})\le Chr(\mathcal{M}_{i})\le k$, so by Theorem \ref{TH_CooperativeFinite} we can choose a function $\Phi _{A}:A\rightarrow \mathcal{K}_{A}$, with $\Phi _{A}(x)\in L(x)$ for all $x\in A$ and $\Phi^{-1}(i)\in \mathcal{I}(\mathcal{M}_{i,A})\subseteq \mathcal{I}(\mathcal{M}_{i})$ for all $i\in \mathcal{K}_{A}$.

     For any $A\in [S]^{<\omega }$, let $T_{A}=\{ X\in [S]^{<\omega }:A\subseteq X\}$, and let us denote $T_{x}=T_{\{x\}}$. Clearly for any $A$, we have $A\in T_{A}$ and for $A_{1},...,A_{n}\in [S]^{<\omega }$, we have $T_{A_{1}}\cap ...\cap T_{A_{n}}=T_{A_{1}\cup ...\cup A_{n}}$. Since $\mathscr{T}=\{T_{A}:A\in [S]^{<\omega }\}\subseteq \mathcal{P}([S]^{<\omega })$ is a set that is closed under finite intersection and $\emptyset \not\in \mathscr{T}$, 
   there is an ultrafilter $\mathscr{U}\subseteq \mathcal{P}([S]^{<\omega })$, such that $\mathscr{T}\subseteq \mathscr{U}$.
  
  Let $x\in S$ be arbitrary. Then for each $A\in T_{x}$, since $x\in A$, the coloring $\Phi _{A}$ is defined on $x$ and $\Phi _{A}(x)\in L(x)$. For each $i\in L(x)$, let $T_{x,i}=\{A\in T_{x}| \Phi _{A}(x)=i\}$. Then the finite union $\bigcup_{i\in L(x)}{T_{x,i}}=T_{x}\in \mathscr{T}\subseteq \mathscr{U}$, so there is a unique $i\in L(x)$, with $T_{x,i}\in \mathscr{U}$. Define the coloring $\Phi _{\mathscr{U}}:S\rightarrow \bigcup_{x\in S}{L(x)}$, be such that for all $x\in S$, $\Phi _{\mathscr{U}}(x)\in L(x)$ is the unique element with $T_{x,\Phi _{\mathscr{U}}(x)}\in \mathscr{U}$.
  
  We need to show that the properties hold for $\Phi _{\mathscr{U}}$. Property 1. is clear from the definition of $\Phi_{\mathscr{U}}$. For property 2.  suppose for contradiction, that there is an $i\in \mathcal{K}$, such that $\Phi _{\mathscr{U}}^{-1}(i)\not\in \mathcal{I}(\mathcal{M}_{i})$. Then there is a circuit $C\in \mathcal{C}(\mathcal{M}_{i})$, such that $C\subseteq \Phi _{\mathscr{U}}^{-1}(i)$. Then for all $x\in C$, we have $i\in L(x)$ and $T_{x,i}\in \mathscr{U}$. Then the finite intersection $\bigcap_{x\in C}{T_{x,i}}\in \mathscr{U}$, so $\bigcap_{x\in C}{T_{x,i}}\neq \emptyset $. Let $A\in \bigcap_{x\in C}{T_{x,i}}$. Then for all $x\in C$, we have  $A\in T_{x,i}$, so $x\in A$ and $\Phi _{A}(x)=i$. Clearly we have $i\in \mathcal{K}_{A}$, but then $C\subseteq \Phi _{A}^{-1}(i)$, that is a contradiction, since $\Phi _{A}^{-1}(i)\in \mathcal{I}(\mathcal{M}_{i})$. Thus, property 2. also holds for $\Phi _{\mathscr{U}}$.
\end{proof}

For the infinite case, we need some preparation.

\begin{definition}\label{DF_MatroidSupportSet}
    Let $\mathcal{M}=(S,r)$ be a finitary matroid, and $X\subseteq S$. We say that $\mathcal{M}$ is suppoted by $X$ if for any $Y\subseteq S$, $Y\cap X\in \mathcal{I}(\mathcal{M})$ implies $Y\in \mathcal{I}(\mathcal{M})$.
\end{definition}

\begin{lemma}\label{LM_SupportClosed}
    If $\mathcal{M}=(S,r)$ is a loop-free finitary matroid supported by $X\subseteq S$ then all subsets $Z\subseteq S-X$ are closed in $\mathcal{M}$.
\end{lemma}

\begin{proof}
    Let $Z_{0}\in [Z]^{<\omega }$ and $y\in S-Z$ arbitrary. Choose $Z'\subseteq Z_{0}$ that is independent and $|Z'|=r(Z_{0})$. Since $Z'\subseteq Z_{0}\subseteq Z$, we have $Z'\cap X=\emptyset$, thus if $y\in X$, we have $(Z'+y)\cap X=\{y\}\in \mathcal{I}(\mathcal{M})$, since $\mathcal{M}$ is loop-free, if $y\not\in X$, we have $(Z'+y)\cap X=\emptyset \in \mathcal{I}(\mathcal{M})$. Since $\mathcal{M}$ is supported by $X$, we have $Z'+y\in \mathcal{I}(\mathcal{M})$, so $r(Z_{0}+y)\ge r(Z'+y)=|Z'+y|=|Z'|+1=r(Z_{0})+1$. Since $y$ was arbitrary, we have $Z$ is closed in $\mathcal{M}$.
\end{proof}

\begin{lemma}\label{LM_SupportContraction}
    If $\mathcal{M}=(S,r)$ is a finitary matroid supported by $X\subseteq S$. Then for all $Z\subseteq S$ with $X\cap Z=\emptyset $ the contraction matroid $\mathcal{M}'=(S-Z,r')$ is the same as the restricted matroid $\mathcal{M}_{S-Z}$.
\end{lemma}    

\begin{proof}
    We will show that $\mathcal{I}(\mathcal{M}')=\mathcal{I}(\mathcal{M}_{S-Z})$. First let $Y\in \mathcal{I}(\mathcal{M}')$. Then for all $A\in [Y]^{<\omega }$, by lemma \ref{MAT_lm:lemma11} we have , that $|A|=r'(A)\le r(A)\le|A|$, so clearly $r(A)=|A|$, thus $Y\in \mathcal{I}(\mathcal{M}_{S-Z})$. Now suppose that $Y\in \mathcal{I}(\mathcal{M}_{S-Z})=\mathcal{I}(\mathcal{M})\cap \mathcal{P}(S-Z)$, and we would like to show that $Y\in \mathcal{I}(\mathcal{M'})$. By lemma \ref{MAT_lm:lemma12}, we only need to show that for any $Y'\in \mathcal{I}(\mathcal{M}')\cap \mathcal{P}(Z)$, we have $Y\cup Y'\in \mathcal{I}(\mathcal{M})$. Since $X\cap Y'\subseteq X\cap Z=\emptyset $, we have $(Y\cup Y')\cap X=Y\cap X\subseteq Y$, thus $(Y\cup Y')\cap X\in \mathcal{I}(\mathcal{M})$, and since $\mathcal{M}$ is supported by $X$, we must have $Y\cup Y'\in \mathcal{I}(\mathcal{M})$, so we are done.
\end{proof}

\begin{lemma}\label{LM_SupportChromatic}
    If $\mathcal{M}=(S,r)$ is a loop-free finitary matroid supported by $X\subseteq S$, then for the restricted matroid we hve $Chr(\mathcal{M}_{X})=Chr(\mathcal{M})$.
\end{lemma}

\begin{proof}
    For one way $Chr(\mathcal{M}_{X})\le Chr(\mathcal{M})$ is clear, as any proper coloring of $\mathcal{M}$ can be restricted to $X$. For the other way, let $\kappa =Chr(\mathcal{M}_{X})$ and let $\Phi :X\rightarrow \kappa $ be a proper coloring of $\mathcal{M}_{X}$. We define the function $\Phi ^{*}:S\rightarrow \kappa $, with $\Phi ^{*}(x)=\Phi (x)$ for $x\in X$, and $\Phi^{*}(x)=0$ for $x\not\in X$. We need to show that this is a proper coloring of $\mathcal{M}$. For $\gamma \in \kappa ,\gamma\neq 0$, we have ${\Phi ^{*}}^{-1}(\gamma )=\Phi ^{-1}(\gamma )\in \mathcal{I}(\mathcal{M})$. For 0, we have ${\Phi^{*}}^{-1}(0)=\Phi ^{-1}(0)\cup (S-X)$. Then ${\Phi^{*}}^{-1}(0)\cap X=\Phi ^{-1}(0)\in \mathcal{I}(\mathcal{M})$ and $\mathcal{M}$ is supported by $X$, so we also have ${\Phi^{*}}^{-1}(0)\in \mathcal{I}(\mathcal{M})$, hence $\Phi ^{*}$ is a proper coloring, $Chr(\mathcal{M})=\kappa $.
\end{proof}

\begin{definition}\label{DF_FreeExtension}
    Let $\mathcal{M}=(S,r)$ be a finitary matroid and $X\subseteq S$, for any $A\in [S]^{<\omega }$, let $r_{<x>}(A)=r(A\cap X)+|A-X|$ and let $\mathcal{M}_{<X>}=(S,r_{<X>})$.
\end{definition}

\begin{lemma}\label{LM_FreeExtension}
    For a finitary matroid $\mathcal{M}=(S,r)$ and $X\subseteq S$, $\mathcal{M}_{<X>}$ is a finitary matroid, that is supported by $X$. Moreover we have $\mathcal{I}(\mathcal{M}_{<X>})\supseteq \mathcal{I}(\mathcal{M})$, and for the restriction, we have $\mathcal{M}_{<X>,X}=\mathcal{M}_{X}$.
\end{lemma}

\begin{proof}
    First we will show that the properties of Definition \ref{MAT_df:finitary-matroid} hold for $\mathcal{M}_{<X>}=(S,r_{<x>})$. 

    1. $r_{<X>}(\emptyset )=r(\emptyset\cap X)+|\emptyset-X|=r(\emptyset)+|\emptyset|=0+0=0$.

    2. If $A\subseteq B\in [S]^{<\omega }$, we have $A\cap X\subseteq B\cap X$ and $A-X\subseteq B-X$, so $r(A\cap X)\le r(B\cap X)$ and $|A-X|\le |B-X|$, hence $r_{<A>}=r(A\cap X)+|A-X|\le r(B\cap X)+|B-X|=r_{<X>}(B)$.

    3. For any $A\in [S]^{<\omega }$, we have $r_{<X>}(A)=r(A\cap X)+|A-X|\le |A\cap X|+|A-X|=|A|$.

    4. For any $A,B\in [S]^{<\omega }$, we have $r_{<X>}(A)+r_{<X>}(B)=r(A\cap X)+|A-X|+r(B\cap X)+|B-X|\ge r((A\cap X)\cap (B\cap X))+r((A\cap X)\cup (B\cap X))+|(A-X)\cap (B-X)|+|(A-X)\cup (B-X)|=r((A\cap B)\cap X)+r((A\cup B)\cap X)+|(A\cap B)-X|+|(A\cup B)-X|=r_{<X>}(A\cap B)+r_{<X>}(A\cup B)$.

    For the restriction $\mathcal{M}_{<X>,X}=\mathcal{M}_{X}$, for any $A\in [X]^{<\omega }$, we have $r_{<X>}(A)=r(A\cap X)+|A-X|=r(A)+0=r(A)$.

    For the support, let $Y\subseteq S$, be such that $Y\cap X\in \mathcal{I}(\mathcal{M}_{<X>})$. Since $Y\cap X\subseteq X$, by the restriction, we have $Y\cap X=\mathcal{I}(\mathcal{M})$. Let $A\in [Y]^{<\omega }$. Then $A\cap X\in [Y\cap X]^{<\omega }$, so $r(A\cap X)=|A\cap X|$. Then we have $r_{<X>}(A)=r(A\cap X)+|A-X|=|A\cap X|+|A-X|=|A|$. Since this holds for all $A\in [Y]^{<\omega }$, we have $Y\in \mathcal{I}(\mathcal{M}_{<X>})$, thus $\mathcal{M}_{<X>}$ is supported by $X$.

    Finally if $Y\in \mathcal{I}(\mathcal{M})$ holds, then $Y\cap X\in \mathcal{I}(\mathcal{M})$ and $Y\cap X\subseteq X$, so by the restriction, we have $Y\cap X\in \mathcal{I}(\mathcal{M}_{<X>})$, so by the support $Y\in \mathcal{I}(\mathcal{M}_{<X>})$. Hence $\mathcal{I}(\mathcal{M})\subseteq \mathcal{I}(\mathcal{M}_{<X>})$ holds.
    
\end{proof}

Now we turn on to the proof of Theorem \ref{TH_CooperativeFinitaryGeneral} in the case when $\kappa $ is an infinite cardinal. Let $\kappa$ be fixed, and we define the statement $Q^{*}(\lambda )$ for any cardinal $\lambda $.

$Q^{*}(\lambda)$ : If $S$ is a base set with $|S|=\lambda $, $\mathcal{K}$ is a color set $(\mathcal{M}_{i})_{i\in \mathcal{K}}$ is a collection of loop-free finitary matroids on the set $S$ and $L:S\rightarrow \mathcal{P}(\mathcal{K})$ is a listing, such that $Chr(\mathcal{M}_{i})\le \kappa $ for all $i\in \mathcal{K}$ and $|L(x)|\ge \kappa $ for all $x\in S$. Then there is a function $\Phi :S\rightarrow \mathcal{K}$, such that

    1. $\Phi(x)\in L(x)$ for all $x\in S$

    2. $\Phi ^{-1}(i)\in \mathcal{I}(\mathcal{M}_{i})$ for all $i\in \mathcal{K}$. 

Clearly Theorem \ref{TH_CooperativeFinitaryGeneral} holds if $Q^{*}(\lambda )$ holds for all cardinals $\lambda $.  We will prove $Q^{*}(\lambda )$ by induction on $\lambda $.

First we will prove $Q^{*}(\lambda )$ for $\lambda \le \kappa $. Let $\lambda \le \kappa $ and $S,\mathcal{K},(\mathcal{M}_{i})_{i\in \mathcal{K}}$ and $L$, be as required, with $|S|=\lambda$. Let $\prec $ be a well-ordering of $S$ in order type $\lambda \le \kappa $. We define $\Phi $ by transfinite induction, since for all $x\in S$, we have $|\{\Phi (y):y\prec x\}|<\lambda \le \kappa \le |L(x)|$, we can choose $\Phi (x)\in L(x)\setminus \{\Phi (y):y\prec x\}$. For this $\Phi $, Property 1. clearly holds, moreover $\Phi $ is one-to-one. Hence for any $i\in \mathcal{K}$, we have either $\Phi ^{-1}(i)=\emptyset$ or $|\Phi^{-1}(i)|=1$, and since $\mathcal{M}_{i}$ is loop-free, we have $\Phi^{-1}(i)\in \mathcal{I}(M_{i})$, so Property 2. also holds.

Now suppose that $\lambda >\kappa $ and $Q^{*}(\mu )$ hold for all $\mu <\lambda $, and we need to prove $Q^{*}(\lambda )$. Let $S$ be a base set with $|S|=\lambda $, $\mathcal{K}$ be a color set and $(\mathcal{M}_{i})_{i\in \mathcal{K}}$ be a collection of loop-free finitary matroids with $Chr(\mathcal{M}_{i})\le \kappa $ for each $i\in \mathcal{K}$, and $L:S\rightarrow \mathcal{P}(\mathcal{K})$ be a listing, with $|L(x)|\ge \kappa $ for $x\in S$. For $i\in \mathcal{K}$, let us denote $L^{-1}[i]=\{x\in S|i\in L(x)\}$. We may have the following assumptions:

a) $|L(x)|=\kappa $ for all $x\in S$

b) $\bigcup_{x\in S}{L(x)}=\mathcal{K}$

c) $\mathcal{M}_{i}$ is supported by $L^{-1}[i]$ for all $i\in \mathcal{K}$.

For a) we can replace each $L(x)$ with $L'(x)\subseteq L(x),|L'(x)|=\kappa $, as in the solution, we would have $\Phi (x)\in L'(x)\subseteq L(x)$. For b) we can replace $\mathcal{K}$ with $\mathcal{K'}=\bigcup_{x\in S}{L(x)}$, as in a solution for $i\in \mathcal{K}\setminus \mathcal{K'}$, we must have $\Phi^{-1}(i)=\emptyset$. For c) replace each matroid $\mathcal{M}_{i}$ with $\mathcal{M}_{i,<L^{-1}[i]>}$. As the $\mathcal{M}_{i}$ matroids are loop-free, we have $\{x\}\in \mathcal{I}(\mathcal{M}_{i})\subseteq \mathcal{I} (\mathcal{M}_{i,<L^{-1}[i]>})$, so the new matroids are also loop-free, and as $\mathcal{M}_{i,<L^{-1}[i]>}$ is supported by $L^{-1}[i]$, we have $Chr(\mathcal{M}_{i,<L^{-1}[i]>})=Chr(\mathcal{M}_{i,<L^{-1}[i]>,L^{-1}[i]})=Chr(\mathcal{M}_{i,L^{-1}[i]})\le Chr(\mathcal{M}_{i})=\kappa $, so the conditions hold. In the solution, since for any $i\in \mathcal{K}$ and $x\in \Phi^{-1}(i)$, we have $\Phi (x)=i$, so $i\in L(x)$, thus $x\in L^{-1}[i]$. This way we have $\Phi^{-1}(i)\subseteq L^{-1}[i]$, so by the restriction, since $\Phi^{-1}(i)\in \mathcal{I}(\mathcal{M}_{i,<L^{-1}[i]>})$, we have $\Phi ^{-1}(i)\in \mathcal{I}(\mathcal{M}_{i})$.

Now under these assumptions, we define some functions. Let $r:\mathcal{K}\times [S]^{<\omega }\rightarrow \omega $ be common rank function defined by $r(i,A)=r_{\mathcal{M}_{i}}(A)$ for $i\in \mathcal{K},A\in [S]^{<\omega }$. Since for all $i\in \mathcal{K} $, we have $Chr(\mathcal{M_{i}})\le \kappa $, let $\Psi _{i}:S\rightarrow \kappa $ be a proper coloring for $\mathcal{M}_{i}$ and let $\Psi: \mathcal{K}\times S\rightarrow \kappa$ be defined as $\Psi (i,x)=\Psi_{i}(x)$. For all $x\in S$, we have $|L(x)|=\kappa $, let $f:S\times \kappa \rightarrow \mathcal{K}$ be a listing function, such as $L(x)=\{f(x,\alpha ):\alpha \in \kappa \}$ for all $x\in S$. Finally by Lemma \ref{MAT_lm:lemma19}, for all $i\in\mathcal{K}$ since $Chr(\mathcal{M}_{i})\le \kappa $, for any $A\in [S]^{<\omega }$, we have that $|\{x\in S|r_{\mathcal{M}_{i}}(A+x)=r_{\mathcal{M}_{i}}(A)\}|\le \kappa $, so we can define a bookkeeping function $h_{i}:[S]^{<\omega }\times \kappa \rightarrow S$, such that for all $A\in [S]^{<\omega }$, we have $\{x\in S|r_{\mathcal{M}_{i}}(A+x)=r_{\mathcal{M}_{i}}(A)\}\subseteq \{h_{i}(A,\alpha ): \alpha \in \kappa\}$. Let $h:\mathcal{K}\times [S]^{<\omega }\times \kappa \rightarrow S$, de defined as $h(i,A,\alpha )=h_{i}(A,\alpha )$.

Now we use elementary submodels similar to section \ref{MAT_sc:inf}. Let $\theta $ be a sufficiantly large regular cardinal. Let $(S_{\alpha })_{\alpha <cf(\lambda )}$ be subsets of $S$, such that $|S_{\alpha }|<\lambda $ for all $\alpha <cf(\lambda )$ and $\bigcup_{\alpha<cf(\lambda )}{S_{\alpha }}=S$. Now we construct elementary submodels $N_{\alpha }\prec H(\theta )$ by induction in the following way. Let $N_{0}\prec H(\theta )$ be such that $\kappa \cup\{S,\mathcal{K},r,\Psi ,L,f,h\}\subseteq N_{0}$ and $|N_{0}|=\kappa $. For successor ordinal, we define $N_{\alpha +1}\prec H(\theta )$ as $N_{\alpha }\cup S_{\alpha }\subseteq N_{\alpha +1}$ and $|N_{\alpha +1}|=|N_{\alpha }\cup S_{\alpha }|$. For limit ordinal $\alpha $, let $N_{\alpha }=\bigcup_{\beta<\alpha }{N_{\beta }}$. This way, clearly for all $\beta<\alpha <cf(\lambda )$, we have $N_{\beta }\subseteq N_{\alpha }$ and $|N_{\alpha }|<\lambda $. We also have that $S=\bigcup_{\alpha <cf(\lambda )\in}{S_{\alpha }}\subseteq \bigcup_{\alpha <cf(\lambda )\in}{N_{\alpha +1}}=\bigcup_{\alpha <cf(\lambda )\in}{N_{\alpha }}$.

For each $\alpha $, let $Z_{\alpha }=S\cap M_{\alpha}$. Then clearly, we have $Z_{\beta }\subseteq Z_{\alpha}$ for $\beta <\alpha $ and $S=\bigcup_{\alpha <cf(\lambda )}{Z_{\alpha }}$. For each $\alpha $, let $\mathcal{K}_{\alpha }=\bigcup _{x\in Z_{\alpha}}{L(x)}$.

\begin{lemma}\label{LM_SubmodelColorSet}
    $\mathcal{K}_{\alpha}\subseteq N_{\alpha}$ for all $\alpha <cf(\lambda )$.
\end{lemma}

\begin{proof}
    Let $i\in \mathcal{K}_{\alpha }$. Then there is some $x\in Z_{\alpha }$ with $i\in L(x)$, so there is some $\gamma \in \kappa $, with $f(x,\gamma )=i$. Since $N_{\alpha }$ is an elementary submodel and $f\in N_{0}\subseteq N_{\alpha }$, $\gamma \in \kappa \subseteq N_{0}\subseteq N_{\alpha}$ and $x\in Z_{\alpha }\subseteq N_{\alpha }$, we have $i=f(x,\gamma )\in N_{\alpha}$.
\end{proof}

\begin{lemma}\label{LM_SubmodelClosed}
    The set $Z_{\alpha }$ is closed in $\mathcal{M}_{i}$ for all $i\in \mathcal{K}$ and $\alpha <cf(\lambda)$.
\end{lemma}

\begin{proof}
    First we look at the case when $i\in \mathcal{K}_{\alpha }$. Suppose for contradiction, that for some $A\in [Z_{\alpha }]^{<\omega }$ and $x\in S-Z_{\alpha }$, we have $r_{\mathcal{M}_{i}}(A+x)=r_{\mathcal{M}_{i}}(A)$. Then there is some $\gamma \in \kappa $, such that $h(i,A,\gamma )=h_{i}(A,\gamma )=x$. Then since $N_{\alpha }$ is an elementary submodel, $h\in N_{0}\subseteq N_{\alpha }$, by Lemma \ref{LM_SubmodelColorSet} $i\in \mathcal{K}_{\alpha }\subseteq N_{\alpha }$, $A$ is a finite subset os $N_{\alpha }$, so it is also in $N_{\alpha }$ and $\gamma \in \kappa \subseteq N_{0}\subseteq N_{\alpha }$, thus $x=h(i,A,\gamma )\in N_{\alpha }$. But since $x\in S$, that would mean $x\in Z_{\alpha }$ in contradiction with $x\in S-Z_{\alpha }$.

    Now suppose $i\not\in \mathcal{K}_{\alpha }$. Then for all $x\in Z_{\alpha }$, we have $i\not\in L(x)$, so $x\not\in L^{-1}[i]$, hence $Z_{\alpha }\cap L^{-1}[i]=\emptyset$. Since $\mathcal{M}_{i}$ is supported by $L^{-1}[i]$, by Lemma \ref{LM_SupportClosed}, $Z_{\alpha }$ is closed in $\mathcal{M}_{i}$.
\end{proof}

Now we need to construct the function $\Phi :S\rightarrow \mathcal{K}$ with the given properties. By induction on $\alpha <cf(\lambda)$, we construct $\Phi _{\alpha }:Z_{\alpha }\rightarrow \mathcal{K}$, such that the following 3 properties hold:

 1. $\Phi_{\alpha }(x)\in L(x)$ for all $x\in Z_{\alpha }$

 2. $\Phi _{\alpha} ^{-1}(i)\in \mathcal{I}(\mathcal{M}_{i})$ for all $i\in \mathcal{K}$. 

 3. $\Phi _{\beta }\subseteq \Phi _{\alpha }$ for $\beta\le \alpha $.

 For $\alpha =0$, since $|Z_{0}|\le |M_{0}|\le \kappa $, we have that $Q^{*}(|Z_{0}|)$ holds, so we can choose $\Phi _{0}$, that statisfies the properties 1. and 2., and property 3. is clear at the initial step.

 For limit ordinal $\alpha $, let $\Phi_{\alpha }=\bigcup_{\beta <\alpha }{\Phi_{\beta }}$. Since 3. holds under $\alpha $, this is a function with $\Dom (\Phi_{\alpha})=\bigcup_{\beta <\alpha }{\Dom (\Phi _{\beta })}=\bigcup_{\beta <\alpha }{Z_{\beta }}=Z_{\alpha }$, and property 3. clearly hold for $\alpha $. Property 1. also holds as for any $x\in Z_{\alpha }$, there is $\beta <\alpha $, with $x\in Z_{\beta }$, so $\Phi _{\alpha }(x)=\Phi_{\beta }(x)\in L(x)$. For property 2. suppose for contradiction, that $\Phi _{\alpha }^{-1}(i)\not\in \mathcal{I}(\mathcal{M}_{i})$ for some $i\in \mathcal{K}$. Then there is a $C\in \mathcal{C}(\mathcal{M}_{i})$ with $C\subseteq \Phi_{\alpha }^{-1}(i)$. Let $C=\{x_{1},...,x_{n}\}$, and for all $1\le j\le n$, choose $\beta _{j}<\alpha $, such that $x_{j}\in Z_{\beta _{j}}$, and let $\beta =\max_{1\le j\le n}{(\beta _{j}})<\alpha $. Then for all $j$, we have $x_{j}\in Z_{\beta _{j}}\subseteq Z_{\beta }$ and $\Phi _{\beta }(x_{j})=i$. But then $C\subseteq \Phi _{\beta }^{-1}(i)$ in contradiction with property 2. for $\beta $ and $i$, thus property 2. most hold for $\alpha $ too.

 Now we turn into the successor step. Suppose, we have constructed $\Phi _{\alpha }$ for some $\alpha <cf(\lambda )$ and we make it for $\alpha +1$. For all $i\in \mathcal{K}$, let $\mathcal{M}_{i}'=\bigslant{\mathcal{M}_{i,Z_{\alpha +1}}}{Z_{\alpha }}$, the matroid, we get from $\mathcal{M}_{i}$ with restricting to $Z_{\alpha +1}$ and contracting $Z_{\alpha }$. These are matroids on the base set $Z_{\alpha +1}-Z_{\alpha }$. By Lemma \ref{LM_SubmodelClosed} for all $i\in \mathcal{K}$, $Z_{\alpha }$ is closed in $\mathcal{M}_{i}$, so by Lemma \ref{MAT_lm:lemma13}, we have that $\mathcal{M}_{i}'$ is loop-free.

 \begin{lemma}\label{LM_ControctionChromatic}
     $Chr(\mathcal{M}_{i}')\le \kappa$ for all $i\in \mathcal{K}$.
 \end{lemma}

 \begin{proof}
     First we look at the case, when $i\not\in \mathcal{K}_{\alpha }$. Then for all $x\in Z_{\alpha }$, we have $i\not\in L(x)$, so $x\not\in L^{-1}[i]$, thus $Z_{\alpha }\cap L^{-1}[i]=\emptyset $. By our assumption $\mathcal{M}_{i}$ is supported by $L^{-1}[i]$, so by Lemma \ref{LM_SupportContraction} , for the contraction, we have $\bigslant{\mathcal{M}_{i}}{Z_{\alpha }}=\mathcal{M}_{i,S-Z_{\alpha }}$, so $\mathcal{M}_{i}'=\mathcal{M}_{i,Z_{\alpha +1}-Z_{\alpha }}$. Hence $Chr(\mathcal{M}_{i}')\le Chr(\mathcal{M}_{i})\le \kappa $.

     Now suppose, that $i\in \mathcal{K}_{\alpha }$, so by Lemma \ref{LM_SubmodelColorSet}, we have $i\in N_{\alpha }$. Then in the model $N_{\alpha }$ for any finite set $C$, we can define $C\in \mathcal{C}(\mathcal{M}_{i})$ in the following way:
     $$C\in \mathcal{C}(\mathcal{M}_{i}):r(i,C)<|C|\wedge (\forall X)((X\subseteq C\wedge X\neq C)\rightarrow r(i,X)=|X|)$$
     This formula has the parameters $i$ and $r$, where $r\in N_{0}\subseteq N_{\alpha }$, so this can be written in $N_{\alpha }$. We will show that the restriction $\Psi_{i}'=\Psi _{i}|_{Z_{\alpha +1}-Z_{\alpha }}:Z_{\alpha +1}-Z_{\alpha }\rightarrow \kappa$ is a proper coloring for $\mathcal{M}_{i}'$.

     Suppose for contradiction, that $\Psi_{i}'^{-1}(\gamma )\not\in \mathcal{I}(\mathcal{M}_{i}')$ for some $\gamma \in \kappa $. Then by Lemma \ref{MAT_lm:lemma12}, there is a $Y\subseteq Z_{\alpha }$ with $Y\in \mathcal{I}(\mathcal{M}_{i})$ and $Y\cup \Psi '^{-1}(\gamma )\not\in \mathcal{I}(\mathcal{M}_{i})$. Let $C\in \mathcal{C}(\mathcal{M}_{i})$ be such that $C\subseteq Y\cup \Psi'^{-1}(\gamma )$. Then, since $C\cap Z_{\alpha }\subseteq Y\in \mathcal{I}(\mathcal{M}_{i})$, we must have $C\not\subseteq Z_{\alpha }$, and for all $x\in C-Z_{\alpha }$, we have $\Psi_{i}(x)=\gamma $. Let
     $$k=\min\{|A|:A\in [Z_{\alpha }]^{<\omega }\wedge \exists C'\in \mathcal{C}(\mathcal{M}_{i}),A\subseteq C'\subseteq Z_{\alpha +1}\wedge C\not\subseteq Z_{\alpha }\wedge \forall x\in C'-A,\Psi_{i}(x)=\gamma\}$$
     This number is well defined, as $A_{0}=C-\Psi_{i}^{-1}(\gamma )$ fits in the definition with $C$. It is also clear, that $k>0$ as in the case $k=0$, we had $C'\subseteq \Psi_{i}^{-1}(\gamma )$ in contradiction with $\Psi _{i}$ is a proper coloring of the matroid $\mathcal{M}_{i}$.

     Let $A\in [Z_{\alpha }]^{<\omega }$ and $C_{1}\in \mathcal{C}(\mathcal{M}_{i})$ be such that $|A|=k$, $A\subseteq C_{1}\subseteq Z_{\alpha +1}$, $C_{1}\not\subseteq Z_{\alpha }$ and $C_{1}-A\subseteq \Psi_{i}^{-1}(\gamma )$. Let $l=|C_{1}-A|=|C_{1}|-k$. Then we can write:
     $$H(\theta )|=\exists x_{1},...,\exists x_{l},\Psi (i,x_{1})=\gamma \wedge ...\wedge \Psi (i,x_{l})=\gamma \wedge A\cup \{x_{1},...,x_{l}\}\in \mathcal{C}(\mathcal{M}_{i})$$
     For the parameters in this expression, we have $i\in N_{\alpha }$ and the circuits of $\mathcal{M}_{i}$ can be written in $N_{\alpha }$, moreover $A\subseteq N_{\alpha }$ is a finite subset, so $A\in N_{\alpha }$,  $\Psi \in N_{0}\subseteq N_{\alpha }$ and $\gamma \in \kappa \subseteq N_{0}\subseteq N_{\alpha }$. Since $N_{\alpha }\prec H(\theta )$ is an elementary submodel, we can write:
     $$N_{\alpha }|=\exists x_{1},...,\exists x_{l},\Psi (i,x_{1})=\gamma \wedge ...\wedge \Psi (i,x_{l})=\gamma \wedge A\cup \{x_{1},...,x_{l}\}\in \mathcal{C}(\mathcal{M}_{i})$$
     Then we have a $C_{2}\in \mathcal{C}(\mathcal{M}_{i})$ with $A\subseteq C_{2}\subseteq (N_{\alpha }\cap S)=Z_{\alpha }$ and $C_{2}-A\subseteq \Psi _{i}^{-1}(\gamma )$. Pick $e\in A\subseteq C_{1}\cap C_{2}$ and $e_{1}\in C_{1}-Z_{\alpha }\subseteq C_{1}-C_{2}$. By Lemma \ref{MAT_lm:lemma2} b), there is a $C_{3}\in \mathcal{C}(\mathcal{M}_{i})$, such that $C_{3}\subseteq (C_{1}\cup C_{2})-e$ and $e_{1}\in C_{3}$. Then we have $C_{3}\subseteq Z_{\alpha +1 }$ and $C_{3}\not\subseteq Z_{\alpha }$. Let $A'=A\cap C_{3}$. Then $A'\subseteq C_{3}$ and for all $x\in C_{3}-A'=C_{3}-A\subseteq (C_{1}\cup C_{2})-A=(C_{1}-A)\cup (C_{2}-A)$, we have $\Psi_{i}(x)=\gamma $. On the other hand, we have $A'\subseteq A-e$, so $|A'|\le |A-e|=k-1$ in contradiction with the definition of $k$. Thus $\Psi _{i}'$ must be a proper coloring for $\mathcal{M}_{i}'$.
 \end{proof}

 Since $|Z_{\alpha +1}-Z_{\alpha }|\le |Z_{\alpha +1}|<\lambda $, we have that $Q^{*}(|Z_{\alpha +1}-Z_{\alpha }|)$ holds, and all of its conditions for the $\mathcal{M}_{i}'$-s and $L$ holds, so there is a function $\Phi _{\alpha }':Z_{\alpha +1}-Z_{\alpha }\rightarrow \mathcal{K}$, such that:

 1. $\Phi_{\alpha }'(x)\in L(x)$ for all $x\in Z_{\alpha +1}-Z_{\alpha }$

 2. $\Phi _{\alpha} '^{-1}(i)\in \mathcal{I}(\mathcal{M}_{i}')$ for all $i\in \mathcal{K}$. 

Let $\Phi _{\alpha +1}=\Phi _{\alpha }\cup \Phi_{\alpha }'$. This is a function with $\Dom(\Phi _{\alpha +1})=\Dom (\Phi _{\alpha } )\cup \Dom (\Phi _{\alpha }')=Z_{\alpha }\cup (Z_{\alpha +1}-Z_{\alpha })=Z_{\alpha +1}$. Property 1. holds, since for all $x\in Z_{\alpha +1 }$, if $x\in Z_{\alpha }$, we have $\Phi _{\alpha +1}(x)=\Phi _{\alpha }(x)\in L(x)$, if $x\in Z_{\alpha +1}-Z_{\alpha }$, we have $\Phi _{\alpha +1}(x)=\Phi _{\alpha }'(x)\in L(x)$. For property 2. let $i\in \mathcal{K}$. Since $\Phi _{\alpha +1}^{-1}(i)=\Phi_{\alpha }^{-1}(i)\cup \Phi _{\alpha }'^{-1}(i)$, by the induction, we have $\Phi_{\alpha }^{-1}(i)\in \mathcal{I}(\mathcal{M}_{i})$ and $\Phi_{\alpha }^{-1}(i)\subseteq Z_{\alpha }$, and in the contraction $\Phi_{\alpha }'^{-1}(i)\in \mathcal{I}(\mathcal{M}_{i}')$, by Lemma \ref{MAT_lm:lemma12}, we have that $\Phi _{\alpha +1}^{-1}(i)\in \mathcal{I}(\mathcal{M}_{i})$. Property 3. is clear, as $\Phi _{\beta }\subseteq \Phi _{\alpha }\subseteq \Phi _{\alpha +1 }$ for all $\beta \le \alpha $. 

Finally, let $\Phi =\bigcup _{\alpha <cf(\lambda )}{\Phi _{\alpha }}$. Then similarly to the limit step, we have that $\Phi $ is a function $\Dom (\Phi )=\bigcup_{\alpha <cf(\lambda )}{\Dom (\Phi _{\alpha })}=\bigcup_{\alpha <cf(\lambda )}{Z_{\alpha }}=S$, and properties 1. and 2. hold for $\Phi $. Then we have done the induction step for $Q^{*}(\lambda )$, proving theorem \ref{TH_CooperativeFinitaryGeneral}.

\end{proof}

\begin{corollary}\label{CO_CooperativeFinitaryGeneral}

 Let $S$ be a base set, $\mathcal{K}$ be a color set $(\mathcal{M}_{i})_{i\in \mathcal{K}}$ be a collection of loop-free finitary matroids on the set $S$ and $L:S\rightarrow \mathcal{P}(\mathcal{K})$ be a listing. Suppose that there is a (finite or infinite) cardinal $\kappa $, such that $Chr(\mathcal{M}_{i,L^{-1}[i]})\le \kappa $ for all $i\in \mathcal{K}$ and $|L(x)|\ge \kappa $ for all $x\in S$. Then there is a function $\Phi :S\rightarrow \mathcal{K}$, such that

    1. $\Phi(x)\in L(x)$ for all $x\in S$

    2. $\Phi ^{-1}(i)\in \mathcal{I}(\mathcal{M}_{i})$ for all $i\in \mathcal{K}$.
    
\end{corollary}

\begin{proof}
   We can replace each matroid $\mathcal{M}_{i}$ with $\mathcal{M}_{i,<L^{-1}[i]>}$. As the $\mathcal{M}_{i}$ matroids are loop-free, we have $\{x\}\in \mathcal{I}(\mathcal{M}_{i})\subseteq \mathcal{I} (\mathcal{M}_{i,<L^{-1}[i]>})$, so the new matroids are also loop-free, and as $\mathcal{M}_{i,<L^{-1}[i]>}$ is supported by $L^{-1}[i]$, by Lemma \ref{LM_SupportChromatic}, we have $Chr(\mathcal{M}_{i,<L^{-1}[i]>})= Chr(\mathcal{M}_{i,L^{-1}[i]})\le \kappa $, so the conditions of Theorem \ref{TH_CooperativeFinitaryGeneral} hold. Applying the theorem, in the solution, since for any $i\in \mathcal{K}$ and $x\in \Phi^{-1}(i)$, we have $\Phi (x)=i$, so $i\in L(x)$, thus $x\in L^{-1}[i]$. This way we have $\Phi^{-1}(i)\subseteq L^{-1}[i]$, so by the restriction, since $\Phi^{-1}(i)\in \mathcal{I}(\mathcal{M}_{i,<L^{-1}[i]>})$, we have $\Phi ^{-1}(i)\in \mathcal{I}(\mathcal{M}_{i})$, thus $\Phi $ meets all criteria.
\end{proof}

\begin{corollary}
     Let $S$ be a base set, $\mathcal{M}$ be a loop-free finitary matroid on $S$, $\mathcal{K}$ be a color set and $L:S\rightarrow \mathcal{P}(\mathcal{K})$ be a listing. Suppose that there is a (finite or infinite) cardinal $\kappa $, such that $Chr(\mathcal{M}_{L^{-1}[i]})\le \kappa $ for all $i\in \mathcal{K}$ and $|L(x)|\ge \kappa $ for all $x\in S$. Then there is a function $\Phi :S\rightarrow \mathcal{K}$, such that

    1. $\Phi(x)\in L(x)$ for all $x\in S$

    2. $\Phi ^{-1}(i)\in \mathcal{I}(\mathcal{M})$ for all $i\in \mathcal{K}$.
\end{corollary}

\section{Partition reduction}\label{MAT_sc:red}

In this section, introduce partitional matroids and a  technique called {\em partition reduction}, that can be used to show relations between chromatic and list chromatic number of finitary matroids. Later we will generalize this concept for matroids that are not necessarily finitary.

\begin{definition}\label{DF_Partition Matroid}
 For a set $S$ a {\em partition} of $S$ is a set $\mathcal{P}$, such that $\emptyset\neq T\subseteq S$ holds for all $T\in \mathcal{P}$, $T\cap T'=\emptyset $ holds for $T\neq T'\in \mathcal{P}$ and $\bigcup{\mathcal{P}}=S$. For $A\in [S]^{<\omega }$, we denote $\mathcal{P}_{A}=\{T\in \mathcal{P}|T\cap 
 A\neq \emptyset\}$ and $r_{\mathcal{P}}(A)=|\mathcal{P}_{A}|$. The pair $\mathcal{M}_{\mathcal{P}}=(S,r_{\mathcal{P}})$ is called the {\em partition matroid } by the partition $\mathcal{P}$.   
\end{definition}

\begin{lemma}\label{LM_PartitionMatroidBasic}
    If $\mathcal{P}$ is a partition of $S$, than $\mathcal{M}_{\mathcal{P}}$ is a finitary matroid over $S$.
\end{lemma}

\begin{proof}
    We need to prove the properties of Definition \ref{MAT_df:finitary-matroid} for $\mathcal{M}_{\mathcal{P}}$

    1. Since for any $T\in \mathcal{P}$, we have $T\cap \emptyset=\emptyset$, we must have $\mathcal{P}_{\emptyset}=\emptyset$, so $r_{\mathcal{P}}(\emptyset)=|\mathcal{P}_{\emptyset}|=|\emptyset |=0$.

    2. Let $A\subseteq B\in [S]^{<\omega }$. Then for any $T\in \mathcal{P}_{A}$, we have $T\cap A\neq \emptyset$, so clearly $T\cap B\neq \emptyset $. Thus $\mathcal{P}_{A}\subseteq \mathcal{P}_{B}$, so $r_{\mathcal{P}}(A)=|\mathcal{P}_{A}|\le |\mathcal{P}_{B}|=r_{\mathcal{P}}(B)$.

    3. Let $A\in [S]^{<\omega }$. Let $f:\mathcal{P}_{A}\rightarrow A$ be a function, such as $f(T)\in T\cap A$ for all $T\in \mathcal{P}_{A}$. Since for $T\neq T'\in \mathcal{P}_{A}$, we have $T\cap T'=\emptyset $, we must have that $f(T)\neq f(T')$. Thus the function $f$ is one-to-one, so $r_{\mathcal{P}}(A)=|\mathcal{P}_{A}|\le |A|$.

    4. Let $A,B\in [S]^{<\omega }$. For any $T\in \mathcal{P}$, we have $T\cap (A\cup B)\neq \emptyset$ if and only if either $T\cap A\neq \emptyset $ or $T\cap B\neq \emptyset $, so $\mathcal{P}_{A\cup B}=\mathcal{P}_{A}\cup \mathcal{P}_{B}$. For the intersection if $T\cap (A\cap B)\neq \emptyset$, then we have $T\cap A\neq \emptyset $ and $T\cap B\neq \emptyset $ (the converse is not necessarily true), so we have $\mathcal{P}_{A\cap B}\subseteq \mathcal{P}_{A}\cap \mathcal{P}_{B}$. Then we have $r_{\mathcal{P}}(A)+r_{\mathcal{P}}(B)=|\mathcal{P}_{A}|+|\mathcal{P}_{B}|=|\mathcal{P}_{A}\cap \mathcal{P}_{B}|+|\mathcal{P}_{A}\cup \mathcal{P}_{B}|\ge |\mathcal{P}_{A\cap B}|+|\mathcal{P}_{A\cup B}|=r_{\mathcal{P}}(A\cap B)+r_{\mathcal{P}}(A\cup B) $.
\end{proof}

\begin{lemma}\label{LM_PartitionMatroidIndependence}
     If $\mathcal{P}$ is a partition of $S$, then for all $X\subseteq S$, we have $X\in \mathcal{I}(\mathcal{M}_{\mathcal{P}})$ if and only if $|T\cap X|\le 1 $ for all $T\in \mathcal{P}$. 
\end{lemma}

\begin{proof}
    First suppose, that $|T\cap X|\le 1 $ holds for all $T\in \mathcal{P}$, and let $A\in [X]^{<\omega }$. For all $x\in A$, let $T(x)$ be the unique element of $\mathcal{P}$ with $x\in T(x)$. Since, $x\in T(x)\cap A$, we have $T(x)\in \mathcal{P}_{A}$. On the other hand if $x\neq y\in A$, since $|T(x)\cap A|\le |T(x)\cap X|\le 1$, $x$ is its only element, so $y\not\in T(x)\cap A$, thus $T(y)\neq T(x)$. Then the function $T:A\rightarrow \mathcal{P}_{A}$ is one-to-one, so we have $r_{\mathcal{P}}(A)=|\mathcal{P}_{A}|\ge |A|$, so $r_{\mathcal{P}}(A)=|A|$. Since $A$ was an arbitrary finite subset of $X$, we have $X\in \mathcal{I}(\mathcal{M}_{\mathcal{P}})$.

    For the converse, suppose, that the condition does not hold, and there is a $T\in \mathcal{P}$, with $x,y\in T\cap X$. Then $\{x,y\}\in [X]^{<\omega }$ and if $T'\in \mathcal{P},T'\neq T$, we have $T'\cap \{x,y\}\subseteq T'\cap T=\emptyset$, so $|\mathcal{P}_{\{x,y\}}|=\{T\}$. Hence $r_{\mathcal{P}}(\{x,y\})=|\mathcal{P}_{\{x,y\}}|=1<2=|\{x,y\}|$, so $X\neq \mathcal{I}(\mathcal{M}_{\mathcal{P}})$.
\end{proof}

\begin{corollary}\label{CO_PartitionMatroidLoopfree}
    If $\mathcal{P}$ is a partition of $S$, than the partition matroid $\mathcal{M}_{\mathcal{P}}$ is loop-free.
\end{corollary}

\begin{lemma}\label{LM_PartitionMatroidChromatic}
     If $\mathcal{P}$ is a partition of $S$, then $Chr(\mathcal{M}_{\mathcal{P}})=List(\mathcal{M_{\mathcal{P}}})=\sup \{|T|:T\in \mathcal{P}\}$
\end{lemma}

\begin{proof}
    Let $\mu =\sup \{|T|:T\in \mathcal{P}\}$. We will show that $\mu \le Chr(\mathcal{M}_{\mathcal{P}})\le List(\mathcal{M}_{\mathcal{P}})\le \mu $.

    For the first inequality, let $\kappa =Chr(\mathcal{\mathcal{M}_{\mathcal{P}}})$ and $\Phi :S\rightarrow \kappa $ be a proper coloring of $\mathcal{M}_{\mathcal{P}}$. Then for any $T\in \mathcal{P}$ if $x\neq y\in T$, by Lemma \ref{LM_PartitionMatroidIndependence} there is no $\gamma \in \kappa $, such that $x,y\in \Phi ^{-1}(\gamma )$, so we must have $\Phi (x)\neq \Phi (y)$. Then the restricted function $\Phi |_{T}:T\rightarrow \kappa $ is one-to-one, we have $|T|\le \kappa $. Since it holds for any $T\in \mathcal{P}$, it holds for the supremum, so $\mu \le \kappa =Chr(\mathcal{M}_{\mathcal{P}})$

    The second inequality $Chr(\mathcal{M}_{\mathcal{P}})\le List(\mathcal{M}_{\mathcal{P}})$ holds in general for all matroids.

    For the third inequality, let $L$ be a listing on $S$, such that $|L(x)|\ge \mu $ for all $x\in S$. For each $T\in \mathcal{P}$, since $|T|\le \mu $, there is a well-order $\prec_{T}$ on $T$ of order type $\le \mu$. We construct a function $\Phi _{T}$ by transfinite induction of $\prec _{T}$. Since for all $x\in T$, we have $|\{\Phi _{T}(y):y\prec _{T}x\}|<\mu $ and $|L(x)|\ge \mu $, we can choose $\Phi _{T}(x)\in L(x)\setminus \{\Phi _{T}(y):y\prec _{T}x\}$. This way clearly $\Phi _{T}(x)\in L(x)$ for all $x\in T$ and $\Phi _{T}$ is one-to-one. Let $\Phi =\bigcup_{T\in \mathcal{P}}{\Phi _{T}}$. This is a function, with $\Dom(\Phi )=\bigcup{\mathcal{P}}=S$, and for all $x\in S$, there is a $T\in \mathcal{P}$, with $x\in T$, so $\Phi (x)=\Phi _{T}(x)\in L(x)$. We need to show, that $\Phi $ is a proper coloring of $\mathcal{M}_{\mathcal{P}}$. For this, let $i\in \bigcup_{x\in S}{L(x)}$. Then for any $T\in \mathcal{P}$, we have $|T\cap \Phi ^{-1}(i)|=|\Phi _{T}^{-1}(i)|\le 1$, as $\Phi _{T}$ is one-to-one. Then by lemma \ref{LM_PartitionMatroidIndependence} $\Phi ^{-1}(i)\in \mathcal{I}(\mathcal{M}_{\mathcal{P}})$, thus $\Phi $ is a proper coloring.
\end{proof}

Now we can define partition reductions of finitary matroids.

\begin{definition}\label{DF_Partition Reduction}
    For a loop-free finitary matroid $\mathcal{M}$ on the set $S$, a {\em partition reduction} of $\mathcal{M}$ is a partition matroid $\mathcal{M}_{\mathcal{P}}$ on the set $S$, with $\mathcal{I}(\mathcal{M}_{\mathcal{P}})\subseteq \mathcal{I}(\mathcal{M})$.
\end{definition}

If $\mathcal{M}_{\mathcal{P}}$ is a partition reduction of $\mathcal{M}$, all proper colorings $\Phi :S\rightarrow \mathcal{K}$ of $\mathcal{M}_{\mathcal{P}}$ are also porper colorings of $\mathcal{M}$, as for any $i\in \mathcal{K}$, we have $\Phi ^{-1}(i)\in \mathcal{I}(\mathcal{M}_{\mathcal{P}})\subseteq \mathcal{I}(\mathcal{M})$. This implies that $Chr(\mathcal{M})\le Chr(\mathcal{M}_{\mathcal{P}})$ and $List(\mathcal{M})\le List(\mathcal{M}_{\mathcal{P}})$ holds.

\begin{definition}\label{DF_ChromaticallyFaithful}
    A partition reduction $\mathcal{M}_{\mathcal{P}}$ of $\mathcal{M}$ is {\em chromatically faithful} if $Chr(\mathcal{M}_{\mathcal{P}})=Chr(\mathcal{M})$.
\end{definition}
    
\begin{lemma}\label{LM_ChromaticallyFaithful}
    If a loop free finitary matroid $\mathcal{M}$ has a chromatically faithful partition reduction, then $Chr(\mathcal{M})=List(\mathcal{M})$.
\end{lemma}

\begin{proof}
    Let $\mathcal{M}_{P}$ be a chromatically faithful partition reduction of $\mathcal{M}$, then we have $Chr(\mathcal{M})\le List(\mathcal{M})\le List(\mathcal{M}_{\mathcal{P}})=Chr(\mathcal{M}_{\mathcal{P}})=Chr(\mathcal{M})$
\end{proof}

Of course, this lemma gives nothing new, as $Chr(\mathcal{M})=List(\mathcal{M})$ holds for all loop-free finitary matroids, by Theorem \ref{MAT_tm:theorem3}, but in fact, the similar argument was, that we used for the infinite chromatic number case.

\begin{lemma}
    If $\mathcal{M}$ is a finitary matroid with $Chr(\mathcal{M})=\kappa$, where $\kappa $ is an infinite cardinal, then it has a chromatically faithful partition reduction.
\end{lemma}

\begin{proof}
    By Theorem \ref{MAT_tm:theorem2}, there is a well-ordered base $(B,\le)$ of $\mathcal{M}$, such that for all $b\in B$, we have $|\{x\in S|M_{B}(x)=b\}|\le \kappa $, where $M_{B}$ is defined by definition \ref{MAT_df:mb}. For each $b\in B$, let $S_{b}=\{x\in S|M_{B}(x)=b\}$ and $\mathcal{P}=\{S_{b}:b\in B\}$. This is clearly a partition. We need to show that $\mathcal{I}(\mathcal{M}_{\mathcal{P}})\subseteq \mathcal{I}(\mathcal{M})$. Suppose for contradiction, that there is an $X\in \mathcal{I}(\mathcal{M}_{\mathcal{P}})\setminus \mathcal{I}(\mathcal{M})$. Then there is a $C\in \mathcal{C}(\mathcal{M})$ with $C\subseteq X$. By Lemma \ref{MAT_lm:lemma17}, there are $x\neq y\in C$ with $M_{B}(x)=M_{B}(y)$. Let $b=M_{B}(x)=M_{B}(y)$, then $x,y\in S_{b}\cap C\subseteq S_{b}\cap X$, so $|S_{b}\cap X|\le 2$, so by Lemma \ref{LM_PartitionMatroidIndependence} $X\not\in \mathcal{I}(\mathcal{M}_{\mathcal{P}})$, that is a contradiction. For the chromatic number, since for all $b\in B$, we have $|S_{b}|\le \kappa $, by lemma \ref{LM_PartitionMatroidChromatic}, we have $Chr(\mathcal{M}_{\mathcal{P}})\le \kappa =Chr(\mathcal{M})$, so $Chr(\mathcal{M}_{\mathcal{P}})=Chr(\mathcal{M})$, the reduction is chromatically faithful.
\end{proof}

This brings us the question, whether all loop-free finitary matroids have chromatically faithful partition reduction, giving a stronger result than Seymour's theorem. Unfortunately this is not true even for finite matroids, as we will see in this example:

\begin{example}\label{EX_K4}
    The graphical matroid of the complete graph $K_{4}$ does not have a chromatically faithful partition reduction.
\end{example}

\begin{proof}
    Let $V$ be the vertex-set with $|V|=4$ and $E$ be the edge set $|E|=6$. Let $\mathcal{M}$ be the graphical matroid of the graph, where each set is independent if it does not contain any circuits. Since we can split $E=E_{1}\cup E_{2}$ to disjoint sets, such that both $E_{1}$ and $E_{2}$ is a path through all vertices, we have $Chr(\mathcal{M})=2$.

    Now suppose, that there is a chromatically faithful partition reduction $\mathcal{M}_{\mathcal{P}}$ of $\mathcal{M}$. Then since $Chr(\mathcal{M}_{\mathcal{P}})=2$, by lemma \ref{LM_PartitionMatroidChromatic}, we have $|T|\le 2$ for all $T\in \mathcal{P}$. Then $\mathcal{P}$ may contain pairs and singletons, but if it contains singletons, their number is even, and if we pair them inside $\mathcal{P}$, we get a new partition $\mathcal{P'}$ that is coarser, so $\mathcal{I}(\mathcal{M}_\mathcal{P}')\subseteq \mathcal{I}(\mathcal{M}_\mathcal{P})\subseteq \mathcal{I}(\mathcal{M})$ and we still have $Chr(\mathcal{M}_{\mathcal{P}})=2$. Hence we may suppose, that $\mathcal{P}$ contains only pairs, so $\mathcal{P}=\{T_{1},T_{2},T_{3}\}$ that are pairwise disjoint 2 element sets of $E$. For all $i$, there is at most 1 vertex $v\in V$, such that $T_{i}$ has 2 edges at $v$, and $|V|=4$, so there is some $v\in V$ where edges of $T_{1},T_{2},T_{3}$ meet. Then the  edge set $F$ of the triangle on the vertices of $V- v$ also contains one edge from all $T_{1},T_{2},T_{3}$, so $F\in \mathcal{I}(\mathcal{M}_{P})$, but since $F$ is a triangle, we have $F\not\in \mathcal{I}(\mathcal{M})$ in contradiction, with $\mathcal{M}_{\mathcal{P}}$ is a partition reduction.
\end{proof}

In the next section, we will show another application of partition reduction for matroids that are not finitary.

\section{Chromatic numbers of duals of finitary matroids}\label{MAT_sc:dual}

Finitary matroids can be the generalisation of finite matroids in many concepts, but they lack a core property, there is no dualisation defined on them.  Bruhn and Diestel in \cite{Br13} defined a more general term of infinite matroinds, that contains all finitary matoids and dualisation can be defined in this broader class of matroids.

\begin{definition}\label{DF_InfiniteMatroid General}
    (Definition 1.1. in \cite{Br13})

    A matroid is a pair $\mathcal{M}=(S,\mathcal{I})$, where $\mathcal{I}\subseteq \mathcal{P}(S)$, and the following properties hold:

    1. $\emptyset\in \mathcal{I}$

    2. If $I\subseteq I'\in \mathcal{I}$, then $I\in \mathcal{I}$

    3. If $I\in \mathcal{I}-\mathcal{I}^{max}$ and $I'\in \mathcal{I}^{max}$, then there is an $x\in I'-I$, such that $I+x\in \mathcal{I}$

    4. For any $I\in \mathcal{I}$ and $X\subseteq S$, the set $\{I'\in \mathcal{I}:I\subseteq I'\subseteq X\}$ has a maximal element by inclusion

    Here $\mathcal{I}^{max}$ denotes the maximal elemnts of $\mathcal{I}$ by inclusion. The property 3. can be replaced by the following:

    3'. For any $I\in \mathcal{I}$ and $I'\in \mathcal{I}^{max}$, there is a $B\in \mathcal{I}^{max}$, such that $I\subseteq B\subseteq I\cup I'$.
\end{definition}

In this definition, the independent sets are $\mathcal{I} (\mathcal{M})=\mathcal{I}$, and the bases are $\mathcal{B}(\mathcal{M})=\mathcal{I}^{max}$. It is clear by property 4. that for any $I\in \mathcal{I}(\mathcal{M})$, there is a $B\in \mathcal{B}(\mathcal{M})$ with $I\subseteq B$. 

By Theorem 3.1. in \cite{Br13} , we can define dualisations of matroids, such as for a matroid $\mathcal{M}$, $\mathcal{M}^{*}$ is the unique matroid with $\mathcal{B}(\mathcal{M}^{*})=\{S-B:B\in \mathcal{B}(\mathcal{M})\}$. By Theorem 3.3. in \cite{Br13} it can be shown, that restriction and contraction of a matroid is also a matroid, and for any matroid $\mathcal{M}$ on the set $S$ and $Z\subseteq S$, we have $\bigslant{\mathcal{M}}{Z}=(\mathcal{M}^{*}_{S-Z})^{*}$.

The circuits in general matoids are defined as $C\in \mathcal{C}(\mathcal{M})$ if $C\not\in \mathcal{I}(\mathcal{M})$ and $C$ is minimal with this property. By Lemma 3.8. in \cite{Br13}, for all $X\subseteq S,X\not\in \mathcal{I}(\mathcal{M})$, there is a $C\in \mathcal{C}(\mathcal{M})$, such that $C\subseteq X$. Moreover, by Lemma 3.11. in \cite{Br13}, if $C\in \mathcal{C}(\mathcal{M})$ and $D\in \mathcal{C}(\mathcal{M}^{*})$, then $|C\cap D|\neq 1$.

In this context, finitary matroids can be defined as matroids, such that all $C\in \mathcal{C}(\mathcal{M})$ are finite, or equivalently matroids in which for any $X\subseteq S$, we have $X\in \mathcal{I}(\mathcal{M})$ if and only if $A\in \mathcal{I}(\mathcal{M})$ for all $A\in [X]^{<\omega }$.

For general matroids, we can define colorings similarly as finitary matroids i.e. a function $\Phi :S\rightarrow \mathcal{K}$ is a proper coloring of $\mathcal{M}$ if for all $i\in \mathcal{K}$, we have $\Phi ^{-1}(i)\in \mathcal{I}(\mathcal{M})$. This way we can define the chromatic and list chromatic numbers $Chr(\mathcal{M})$ and $List(\mathcal{M})$ for general matroids.

Partitional matroids are itself finitary, but for a general matroid $\mathcal{M}$, we can define that a partitional matroid $\mathcal{M}_{\mathcal{P}}$ is a partition reduction of $\mathcal{M}$ if $\mathcal{I}(\mathcal{M}_{\mathcal{P}})\subseteq \mathcal{I}(\mathcal{M})$, and $\mathcal{M}_{\mathcal{P}}$ is chromatically faithful if $Chr(\mathcal{M}_{\mathcal{P}})=Chr(\mathcal{M})$. Lemma \ref{LM_ChromaticallyFaithful} can also be used for general matroids, that is not void in this case.

For finitary matroids, the chromatic number can be arbitrarily big, however in case of duals of finitary matroids, this must be countable.

\begin{definition}\label{DF_FreeMatroid}
    A matroid $\mathcal{M}$ on set $S$ is {\em free}, if $\mathcal{I}(\mathcal{M})=\mathcal{P}(S)$ or equivalently $\mathcal{C}(\mathcal{M})=\emptyset$.
\end{definition}

\begin{lemma}\label{LM_ContractionNoFree}
    If $\mathcal{M}$ is a finitary matroid on $S$, such that its dual $\mathcal{M}^{*}$ is loop-free, then for any $Z\subseteq S,Z\neq S$, the contraction matroid $\bigslant{\mathcal{M}}{Z}$ is not free.
\end{lemma}

\begin{proof}
    Suppose for contradiction, that $\bigslant{\mathcal{M}}{Z}$ is free for some $Z$, and let $x\in S-Z$. Since $\mathcal{M}^{*}$ is loop-free, we have $\{x\}\in \mathcal{I}(\mathcal{M}^{*})$, so there is a $B^{*}\in \mathcal{B}(\mathcal{M}^{*})$ with $\{x\}\subseteq B^{*}$, so $x\in B^{*}$. Let $B=S-B^{*}\in \mathcal{B}(\mathcal{M})$, then $x\not\in B$. Since $(B+x)\cap Z=B\cap Z\subseteq B\in \mathcal{I}(\mathcal{M})$ and $(B+x)-Z\in \mathcal{I}(\bigslant{\mathcal{M}}{Z})$ as it is a free matroid, by Lemma \ref{MAT_lm:lemma12}, we have that $B+x\in \mathcal{I}(\mathcal{M})$ in contradiction with the maximality of $B$.
\end{proof}

\begin{theorem}\label{TH_FinitaryDualCountable}
    If $\mathcal{M}$ is a finitary matoid on some base set $S$, such that its dual $\mathcal{M}^{*}$ is loop-free, then $Chr(\mathcal{M}^{*})\le List(\mathcal{M}^{*})\le \omega $.
\end{theorem}

\begin{proof}
    By transfinite induction, we define subsets $C_{\alpha }\subseteq S$ for all ordinals $\alpha $. For every $\alpha $ if $\bigcup_{\beta <\alpha }{C_{\beta}}\neq S$, then let $C_{\alpha}\in \mathcal{C}(\bigslant{\mathcal{M}}{\bigcup_{\beta<\alpha }{C_{\beta}}})$ (by Lemma \ref{LM_ContractionNoFree} the contraction is not free, so it contains a circuit), if $\bigcup_{\beta <\alpha }{C_{\beta}}= S$, then let $C_{\alpha }=\emptyset $. Clearly, by this definition, these are subsets of $S$ and pairwise disjoint, so they cannot be all non empty. Let $\rho $ be the smallest ordinal, such that $C_{\rho }=\emptyset $. Then by the construction $\bigcup _{\alpha <\rho}{C_{\alpha }}=S$, and $C_{\alpha }\neq \emptyset $ for $\alpha <\rho$.

    Let $\mathcal{P}=\{C_{\alpha }:\alpha <\rho \}$. Then this is a partition of $S$ and since for all $\alpha $, $C_{\alpha }$ is a circuit in a contraction of a finitary martoid, that is itself finitary, we have that $|C_{\alpha }|<\omega $ for all $\alpha <\rho$. Then by Lemma \ref{LM_PartitionMatroidChromatic}, we have $Chr(\mathcal{M}_{\mathcal{P}})=List(\mathcal{M}_{\mathcal{P}})=sup\{|C_{\alpha }|:\alpha <\rho\}\le \omega $.

    We will show, that $\mathcal{M}_{\mathcal{P}}$ is a partition reduction of $\mathcal{M}^{*}$. Suppose for contradiction, that there is some $X\in \mathcal{I}(\mathcal{M}_{\mathcal{P}})-\mathcal{I}(\mathcal{M}^{*})$ and let $D\in \mathcal{C}(\mathcal{M}^{*})$ be such that $D\subseteq X$. Let $\alpha $ be the smallest ordinal with $C_{\alpha }\cap D\neq \emptyset$. Then $D\subseteq S-\bigcup_{\beta <\alpha }{C_{\beta}}$, so $D\in \mathcal{C}(\mathcal{M}^{*}_{S-\bigcup_{\beta <\alpha }{C_{\beta}}})=\mathcal{C}((\mathcal{M}^{*}_{S-\bigcup_{\beta <\alpha }{C_{\beta}}})^{**})=\mathcal{C}((\bigslant{\mathcal{M}}{\bigcup_{\beta <\alpha }{C_{\beta}}})^{*})$. Since $C_{\alpha }\in \mathcal{C}(\bigslant{\mathcal{M}}{\bigcup_{\beta <\alpha }{C_{\beta}}})$, by Lemma 3.11. in \cite{Br13}, we have that $|C_{\alpha }\cap D|\neq 1$, so $|C_{\alpha }\cap X|\ge |C_{\alpha }\cap D|\ge 2$, but by Lemma \ref{LM_PartitionMatroidIndependence} this is in contradiction with $X\in \mathcal{M}_{\mathcal{P}}$.

    Hence, we have $Chr(\mathcal{M}^{*})\le List(\mathcal{M}^{*})\le List(\mathcal{M}_{\mathcal{P}})\le \omega $.
\end{proof}

\begin{corollary}\label{CO_ChromaticDual}
    If $\mathcal{M}$ is a loop-free matroid, such that $\mathcal{M}$ is a dual of a finitary matroid and $Chr(\mathcal{M})$ is infinite, then $Chr(\mathcal{M})=List(\mathcal{M})$.
\end{corollary}

\section{Problems}

\begin{problem}
    Is $Chr(\mathcal{M})=List(\mathcal{M})$ holds for all loop free infinite matroids?
\end{problem}

There is either a more specific case, that we don't know:

\begin{problem}
   Is $Chr(\mathcal{M})=List(\mathcal{M})$ holds for loop free matroids, that are duals of finitary matroids and $Chr(\mathcal{M})<\omega$? 
\end{problem}

\begin{problem}
    Do all loop-free infinite matroids with infinite chromatic number $Chr(\mathcal{M})$ have a chromatically faithful partition reduction?
\end{problem}

 \bibliographystyle{plain}

\end{document}